\documentclass[10pt]{article}

\usepackage{auto-pst-pdf}
\usepackage{amsmath, latexsym, amsfonts, amssymb, amsthm}
\usepackage{graphicx, color,hyperref,dsfont,epsfig,caption,wrapfig,subfig}
\usepackage{pstricks}
\usepackage{datetime}
\usepackage{movie15}
\hypersetup{
    linktoc=page,
    linkcolor=red,          
    citecolor=blue,        
    filecolor=blue,      
    urlcolor=cyan,
    colorlinks=true           
}

\numberwithin{equation}{section}

\newtheorem{theorem}{Theorem}[section]
\newtheorem{lemma}[theorem]{Lemma}
\newtheorem{proposition}[theorem]{Proposition}
\newtheorem{corollary}[theorem]{Corollary}

\newtheorem{remark}[theorem]{Remark}
\newtheorem{definition}[theorem]{Definition}

\newcounter{conj}
\newtheorem{conjecture}[conj]{Conjecture}

\theoremstyle{definition}

\renewcommand{\tilde}{\widetilde}          
\DeclareMathSymbol{\leqslant}{\mathalpha}{AMSa}{"36} 
\DeclareMathSymbol{\geqslant}{\mathalpha}{AMSa}{"3E} 
\DeclareMathSymbol{\eset}{\mathalpha}{AMSb}{"3F}     
\renewcommand{\leq}{\;\leqslant\;}                   
\renewcommand{\geq}{\;\geqslant\;}                   

\newcommand{\C}{\mathbb{C}}
\newcommand{\R}{\mathbb{R}}

\newcommand{\E}{\mathds{E}}
\renewcommand{\P}{\mathds{P}}

\usepackage{xparse}

\DeclareDocumentCommand \Pmp { m m o} {
\IfNoValueTF{#3}
{P_{#1}^{#2}}
{P_{#1}^{#2}\left(#3\right)}
}

\DeclareDocumentCommand \Emp { m m o} {
\IfNoValueTF{#3}
{E_{#1}^{#2}}
{E_{#1}^{#2}\left[#3\right]}
}

\DeclareDocumentCommand \Pbr { m m m m o } {
\IfNoValueTF{#5}
{P_{#1}^{#2\stackrel{#4}{\rightarrow} #3}}
{P_{#1}^{#2\stackrel{#4}{\rightarrow} #3}\left(#5\right)}
}

\DeclareDocumentCommand \Ebr { m m m m o } {
\IfNoValueTF{#5}
{E_{#1}^{#2\stackrel{#4}{\rightarrow} #3}}
{E_{#1}^{#2\stackrel{#4}{\rightarrow} #3}\left[#5\right]}
}

\def\S{\mathbb{S}}

\def\bi{\begin{itemize}}
\def\ei{\end{itemize}}

\def\bnum{\begin{enumerate}}
\def\enum{\end{enumerate}}
\def\<#1{\langle #1 \rangle}

\title{Liouville Quantum Gravity on the Riemann sphere}

\author{ Fran\c{c}ois David \footnote{Institut de Physique Th\'eorique,
CNRS, URA 2306, CEA, IPhT, Gif-sur-Yvette, France.}, Antti Kupiainen \footnote{University of Helsinki, Department of Mathematics and Statistics, P.O.
Finland. Supported by the Academy of Finland}, R\'emi Rhodes \footnote{Universit{\'e} Paris-Est Marne la Vall\'ee, LAMA, Champs sur Marne, France.} \footnotetext[2]{Partially supported by grant ANR-11-JCJC  CHAMU.},
 Vincent Vargas \footnote{ENS Ulm, DMA, 45 rue d'Ulm,  75005 Paris, France.} }

\begin{document}

\maketitle

\begin{abstract}
In this paper, we rigorously construct Liouville Quantum Field Theory  
on the Riemann sphere introduced in the 1981 seminal work by Polyakov. 
We establish some of its fundamental properties like conformal covariance under PSL$_2(\C)$-action, Seiberg bounds, KPZ scaling laws,  KPZ formula and the Weyl anomaly 
formula. We also make precise conjectures about the relationship of
the theory to scaling limits of random planar maps conformally embedded onto the sphere. 
\end{abstract}

\begin{center}
\end{center}
\footnotesize



\noindent{\bf Key words or phrases:}  Liouville Quantum Gravity, quantum field theory, Gaussian multiplicative chaos, KPZ formula, KPZ scaling laws, Polyakov formula.

\noindent{\bf MSC 2000 subject classifications:  81T40,  81T20, 60D05.}    

\normalsize



\section{Introduction}

The two dimensional  Liouville Quantum Field  Theory  (Liouville QFT for short, or also LQG \footnote{not to be confused with Loop Quantum Gravity, another approach to quantize gravity in 3 and 4 dimensions...} as a short-cut for Liouville Quantum Gravity  as is now usual in the mathematics literature) was introduced by A. Polyakov in 1981 \cite{Pol} as a model for quantizing the bosonic string in the conformal gauge  and gravity in  two space-time dimensions. Liouville QFT  is one of the most important two dimensional Conformal Field Theories (CFT).

{\it Classical} Liouville theory is a theory of Riemannian metrics $g$ on a two dimensional surface $\Sigma$. One considers metrics $g=e^{\gamma X}\hat g$ where $\hat g$ is some
fixed smooth "reference" metric, $X:\Sigma\to\R$ is a deterministic function and $\gamma$ is a real parameter. The Liouville action functional is then defined as
\begin{equation}\label{actionfalse}
S(X,\hat{g}):= \frac{1}{4\pi}\int_{\Sigma}\big(|\partial^{\hat{g}}X |^2+QR_{\hat{g}} X  +4\pi \mu e^{\gamma X  }\big)\,\lambda_{\hat{g}} ,   
\end{equation}
where $\partial^{\hat{g}}$, $R_{\hat{g}} $ and $\lambda_{\hat{g}}$ respectively stand for the gradient, Ricci scalar curvature and volume form in the metric $ \hat{g}$. The parameter $\mu> 0$ is the analog of a``cosmological constant''  in two dimensional gravity and  $Q$ is a real parameter. For the particular value 
$$Q=\frac{2}{\gamma}$$
this action functional  is {\it conformally invariant}. This means that if we choose a complex coordinate $z$ so that the metric is given as $\hat g=dz^2$ then \eqref{actionfalse}
is invariant under the simultaneous change of coordinates $z=f(w)$ and shift in the field
\begin{equation}
\label{classconf}
X= X'\circ f + Q \log |f'| .
\end{equation}
 In that case its
extrema  are given by solutions of the classical  {\it Liouville equation} 
$$R_{e^{\gamma X}\hat{g}} =-2\pi\mu\gamma^2.$$
Such solutions define metrics  $e^{\gamma X}\hat g$ with  constant negative curvature and lead to the uniformisation theorem of Riemann surfaces.

In the {\it quantum} (or {\it probabilistic}) Liouville theory the field $X$ becomes a random field
with law given heuristically  in terms of a  functional integral 
\begin{equation}\label{pathintegral}
\E[F(X)]=Z^{-1}\int F(X)e^{-S(X,\hat{g})}DX
\end{equation}
where $Z$ is a normalization constant and $DX$ stands for a formal uniform measure on 
some space of maps  $X:\Sigma\to\R$.
We stress that for $\mu>0$ this field is non Gaussian whereas for $\mu=0$ it is Gaussian, in which case it is known under the name Gaussian Free Field (GFF), free referring to the fact that the field is "free" of  interactions (here the term $4\pi\mu\, \mathrm{e}^{\gamma X}$). The aim of this paper is to make rigorous sense of the heuristic expression \eqref{pathintegral} and study its properties.


The quantum Liouville theory is a {\it Conformal Field Theory}. 
This means that we expect there to be a sense in which the random field $X$ is invariant under the conformal transformations \eqref{classconf} 
(see section \ref{Conformal covariance}), {\it however} with a renormalized value of the $Q$ parameter
 $$Q=\frac{2}{\gamma} +\frac{\gamma}{2}.$$
 Note that $Q$ is invariant under $\gamma\to 4/\gamma$ and  in the standard branch of Liouville theory, the parameter $\gamma$ belongs to the interval $]0,2]$ (see Section 4 for discussion of $\gamma\geq 2$). 
 Conformal Field Theories are characterized  by the central charge $c\in\R$  that reflects the way the theory reacts to changes of the background metric (see section \ref{sec:Weyl}). For the Liouville quantum theory,  the central charge  is $c=1+6 Q^2$: thus it  can range continuously in the interval $[25,+\infty[$ and this is one of the interesting features of this theory. 

Since its introduction, Liouville QFT has been and is still much studied in theoretical physics, in the context of integrable systems and conformal field theories, of string theories, of quantum gravity, for its relation with random matrix models  and topological gravity (see \cite{nakayama} for a review), and more recently in the context of its relations with 4 dimensional supersymmetric gauge theories and the AGT conjecture \cite{AGT}.

Liouville theory has also raised recently much interest in mathematics and theoretical physics in the (slightly different) context of probability theory and random geometry where the conjectured link between large planar maps and LQG is intensively studied (see Section \ref{sec:maps} for a discussion on this point).  Up to now such studies have exclusively focused on the incarnation of the Liouville theory as a free field theory  where the parameter $\mu$ is set to zero\footnote{In the mathematics literature, one speaks of critical LQG when $\mu=0$ though the terminology is misleading because non critical LQG, which is the object of this work, is also a CFT.}: see for instance \cite{cf:DuSh,Rnew10} and the review \cite{bourbaki}. Within this framework, the Liouville measure $\mathrm{e}^{\gamma X}\lambda_{\hat{g}}$ is formally the exponential of the GFF and is mathematically defined  via Kahane's theory of Gaussian multiplicative chaos   \cite{cf:Kah} for $\gamma\in]0,2[$. It is then possible to study in depth the properties of the measure in relation with SLE curves or geometrical objects in the plane that can be constructed out of the GFF \cite{AJKS,dubedat,zipper}. In particular, in this geometrical and probabilistic context, a precise mathematical formulation of the KPZ scaling relations can be given \cite{Aru, cf:DuSh,Rnew10}. 

Let us also mention that defining a random metric is an important open problem in the field and steps towards this problem have been achieved in  \cite{Curien} and \cite{MS} in the special case $\gamma=\sqrt{8/3}$: recall that for this value of $\gamma$ (and in the context $\mu=0$) the work \cite{MS} constructs a random growth process which is conjectured to be the growth of balls of a metric space formally corresponding to the "tangent plane" of LQG. This metric space is supposed to correspond to the conformal embedding in the plane of the so-called Brownian plane, recently constructed in \cite{CuLe}. Also, as originally suggested in \cite{David-KPZ}, one can define rigorously the associated diffusion process called Liouville Brownian motion \cite{GRV} (see also \cite{berest} for a construction starting from one point). This has led to further understanding of the geometry of ($\mu=0$) LQG via heat kernel techniques, see \cite{andres,heatKPZ,GRV-FD,MRVZ,spectral} for recent progresses.

Treating the Liouville theory as a GFF (thus setting $\mu=0$)
 is justified in some cases. It was used in the physics literature in the original derivation of KPZ exponents 
(see the seminal work \cite{cf:KPZ} and also \cite{cf:Da,DistKa} for the framework considered here) and it is  the basis of many formal calculations  \cite{DiFran,GouLi,GuTriWi} of the correlation functions of the Liouville theory, i.e. expectations of product of {\it vertex operators} which are random fields of the form 
\begin{equation}\label{vertex}
V_\alpha(x)=\mathrm{e}^{\alpha X(x)}
\end{equation}
(properly renormalized, see Section \ref{sec:gmc}). Indeed, if one performs a formal expansion of the interaction term $4 \pi \mu e^{\gamma X}$ in the formula \eqref{pathintegral} in powers of $\mu$ one ends up computing expectation values in the GFF of products of fields \eqref{vertex} integrated over their location which may be
calculated in closed form 
 \cite{DotFateev1985}. 
Such calculations lead in particular to the famous DOZZ formula for the 3-point correlation functions of Liouville theory on the sphere \cite{Do,ZZ} (see section \ref{DOZZ} for further explanations).
Thanks to these calculations, many checks have been done between the results of Liouville theory for the correlation functions and  corresponding calculations using random matrix models and integrable hierarchies \cite{nakayama}.

%

Nevertheless for many questions the "interaction'' exponential term has to be taken into account. This is for instance the case for the open string (Liouville theory in the disk) where the negative curvature metric and the boundary conditions play an essential role. The purpose of this paper is precisely to define the full Liouville theory for all $\mu>0$ in the simple case of the theory defined on the Riemann sphere $\S^2$ (the theory in the disk can be defined along the same lines  but details of the construction will appear elsewhere). 
The small scale properties of  Liouville field theory are relatively simple:  a simple normal ordering renders the interaction term well defined for $\gamma<2$. 
However, the
theory has unconventional properties in large scales due to a neutral direction ("zero mode") in the integral \eqref{pathintegral}.

We will construct the general $k$-point correlation functions of vertex operators \eqref{vertex} on the sphere
satisfying the so called {\it Seiberg bounds}  \cite{seiberg}:
\begin{equation}\label{introS}
\sum_{i=1}^k \alpha_i>2Q\quad \text{ and }\quad \alpha_i<Q, \,\,\, \,\, \forall i.
\end{equation}
We will also study the conformal invariance and $\mu$-dependence of these correlation functions and study the associated Liouville measure.  In particular, we establish the well known {\it KPZ scaling laws} (see \cite{cf:KPZ,nakayama,witten}) on the $\mu$-dependence. Finally we determine the way the correlations behave under conformal changes of metrics, known as 
the {\it Weyl anomaly} formula (see \cite{brown76,CoJa,CaDuff,duff77,Pol} for early references on the scale and Weyl anomalies in the physics literature and  \cite{R-S,sarnak} on related mathematical works)
thereby recovering  
$c_L=1+6Q^2$ as the {central charge} of the Liouville theory. Finally, we discuss possible approaches of the $\gamma\geq 2$ branches of LQG and formulate precise conjectures on the relationship between LQG and scaling limits of planar maps.

 Our results should not appear as a surprise for theoretical physicists as we recover (in a rigorous setting) many known properties of LQG but they are the first rigorous probabilistic results about the full Liouville theory (on the sphere), as it was introduced by Polyakov in his 1981 seminal paper \cite{Pol}.

\section{Background}\label{sec:backgr}

 \subsection{Metrics on  $\R^2$}\label{sphermetric}    
The  sphere $\S^2$ can be mapped by stereographic projection to the plane which we view both as   ${\R^2}$ and as $\C$. Given a Riemannian metric on $\R^2$
we will denote by $\partial^g$ the gradient, $\triangle_g$ the Laplace-Beltrami operator, $R_g=-\triangle_g\ln \sqrt{\det g}$ the Ricci scalar curvature and $\lambda_g$ the volume form in the metric $g$. When no index is given, this means that the object has to be understood in terms of the usual Euclidean metric on the plane (i.e. $\partial$, $\triangle$, $R$ and $\lambda$). 

We take as the  background metric in \eqref{actionfalse} the  spherical metric on $\S^2$ which becomes on $\R^2$ and on $\C$
$$\hat{g}=\frac{4}{(1+|x|^2)^2}dx^2=\frac{2}{(1+\bar zz)^2}(dz\otimes d\bar z+d\bar z\otimes dz).$$
Its Ricci scalar curvature is $R_{\hat g}=2$ (its Gaussian curvature is $1$) and the volume $\int_{\R^2}\lambda_{\hat g}=4\pi$.

We let $\bar C({\R^2})$ stand for the space of continuous functions on $\R^2$ admitting a finite limit at infinity. In the same way, $\bar C^k({\R^2})$ for $k\geq 1$ stands for the space of $k$-times differentiable functions  on $\R^2$ such that all the derivatives up to order $k$ belong to $\bar C({\R^2})$. 
We say a metric $g=g(x)dx^2$ is conformally equivalent to $\hat g$ if 
$$g(x)=e^{\varphi(x)}\hat{g}(x)
$$ 
with $\varphi\in \bar C^2({\R^2})$  such that $\int_{\R^2}|\partial \varphi|^2\,d\lambda<\infty$. 
We often identify the metrics $g$ with their densities $g(x)$ (or $g(z)$) with respect to the Euclidean metric.  
The curvature $R_g$ can be obtained from the curvature relation
\begin{equation}\label{curvature}
R_g=e^{-\varphi}\big(R_{\hat{g}}-\Delta_{\hat{g}}\varphi\big).
\end{equation}
In what follows, we will denote by $m_{g}(h)$ the mean value of $h$ in the metric $g$, that is
\begin{equation}
m_{g}(h)=\frac{1}{\lambda_g(\R^2)}\int_{\R^2}h\,d\lambda_g.
\end{equation}

Given any metric $g$ conformally equivalent to the spherical metric, one can consider the Sobolev space $H^1(\R^2,g)$, which is the closure of $\bar C^\infty({\R^2})$ with respect to the Hilbert-norm
\begin{equation}\label{Hnorm}
\int_{\R^2}h^2\,d\lambda_{g}+\int_{\R^2}|\partial^g h|^2\,d\lambda_g
\end{equation}
Note that the Dirichlet energy is independent on the metric:
\begin{equation}\label{diri}
\int_{\R^2}|\partial^g h|^2\,d\lambda_g=\int_{\R^2}|\partial h|^2\,d\lambda.
\end{equation}

\subsection{Gaussian free fields}\label{sub:free}

 The purpose of this section is to give a precise meaning to the expression \eqref{pathintegral} in the absence of the $\mu$ and  $Q$ terms i.e. we want to give a precise meaning to the measure formally given by 
\begin{equation}\label{actionformal}
\exp\Big(-\frac{1}{4\pi}\int_{\R^2}  |\partial^{g}X|^2 \lambda_{g} \Big)DX
\end{equation}
where  $g$ is any metric conformally equivalent to the spherical one and $DX$ stands for a "uniform measure" on some space of  $X:\R^2\to\R$. Obviously \eqref{actionformal} should be defined in terms of a Gaussian measure. However, there is an important twist in that we want to include  in the integration domain the constant functions for which the exponent in \eqref{actionformal} vanishes.  This means that the resulting measure will not be a probability measure. Before giving the precise mathematical definition, we choose to explain first the motivations for the forthcoming definitions.

{\bf Heuristic explanation.} By \eqref{diri} the density in \eqref{actionformal} is independent on $g$ and so one can recognize it as a formal density for the Gaussian Free Field (GFF)
 i.e. a centered Gaussian field with covariance structure
\begin{equation}
\E[X(x)X(y)]=\ln\frac{1}{|x-y|}
\end{equation}
(for references to such  log-correlated fields  see  see \cite{dubedat,gaw,She07, DRSVreview}).
This field is defined only up to a constant. One  way to fix the constant is to consider its restriction to the space of test functions $f$ with vanishing mean $\int_{\R^2}fd\lambda=0$ (see \cite{DRSVreview}). This is not the approach that we will develop here. Given a metric $g$ conformally equivalent to $\hat g$, we will rather consider a  field $X_g$ conditioned on having vanishing mean in the metric $g$. Heuristically,  we have 
\begin{equation}
X_{g}=X-m_{g}(X).
\end{equation}
The constant has thus been   fixed by imposing the condition
\begin{equation}
\int X_{g} \,d\lambda_{g}=0. 
\end{equation}
Though this description is not rigorous as the field $X$ does not exist as a function, each field $X_g$ is perfectly defined through its   covariance which is explicitly given
\begin{align}\label{covXphi}
G_{g}(x,y)&:=\E[X_{g}(x)X_{g}(y)]\\
=&\ln\frac{1}{|x-y|}-m_{g}(\ln\frac{1}{|x-\cdot|} )-m_{g}(\ln\frac{1}{|y-\cdot|} )+\theta_{g},\nonumber
\end{align}
with
\begin{equation}\label{theta}
\theta_{g}:=\frac{1}{\lambda_{g}(\R^2)^2}\iint_{\R^2\times\R^2}\ln\frac{1}{|z-z'|}\lambda_{g}(dz)\lambda_{g}(dz').
\end{equation}
It is then plain to check that $X_g$ is a Gaussian Free Field with vanishing $\lambda_g$-mean on the sphere.\qed

\bigskip Therefore we introduce the following definition
\begin{definition}
For each metric $g$ conformally equivalent to $\hat{g}$, we consider a Gaussian Free Field $X_g$ with vanishing $\lambda_g$-mean on the sphere,  that is a centered  Gaussian random distribution with covariance kernel given by the Green function $G_g$ of the problem
\begin{equation}\nonumber
\triangle_g u=-2\pi f \quad\text{on }\R^2,\quad  \int_{\R^2} u \,d\lambda_g =0
\end{equation}
i.e.
\begin{equation}
u=\int G_g(\cdot,z)f(z)\lambda_g(dz):=G_gf.
\end{equation}
\end{definition}

By a straightforward adaptation of \cite{dubedat,She07} one can show that $X_g$ lives almost surely   in the dual space $H^{-1}(\R^2,g)$ of $H^{1}(\R^2,g)$,  and this space does not depend on the choice of the metric $g$ in the conformal equivalence class of $\hat{g}$. We state the following classical result on the Green function $G_g$ (see the appendix for a short proof)
\begin{proposition}{\bf (Conformal covariance)}\label{GreenCC}
Let $\psi$ be a M\"obius transform of the sphere  and consider the metric $g_\psi(z)=|\psi'(z)|^2g(\psi(z))$. We have
$$G_{g_\psi}(x,y)=G_g(\psi(x),\psi(y)). $$
\end{proposition} 
Furthermore, a simple check of covariance structure with the help of \eqref{covXphi} entails
\begin{proposition}{\bf (Rule for changing metrics)}\label{changemetric}
Let the metrics $g,g'$ be conformally equivalent to the spherical metric. Then we have the following equality in law
$$X_g-m_{g'}(X_g)\stackrel{law}{=}X_{g'}.$$
\end{proposition} 

Specializing to the round metric, let us register the explicit formula 
\begin{equation}\label{hatGformula}
G_{\hat g}(z,z')=\ln\frac{1}{|z-z'|}-\frac{1}{4}(\ln\hat g(z)+\ln\hat g(z'))-\frac{1}{2}
\end{equation}
and the transformation rule under M\"obius maps
\begin{equation}\label{Grule}
G_{\hat g}(\psi(z),\psi(z'))=G_{e^\phi\hat g}(z,z')=G_{\hat g}(z,z')-\frac{1}{4}(\phi(z)+\phi(z'))
\end{equation}
where $e^\phi=\hat g_\psi/\hat g$ (see Appendix).

All these GFFs $X_g$ ($g$ conformally equivalent to $\hat{g}$) may be thought of as  centerings in $\lambda_g$-mean  of the same field. They all differ by a constant. To absorb the dependence on the constant, we tensorize the law $\P_g$ of the field $X_g$ 
with the Lebesgue measure $dc$ on $\R$ and we consider the image of the measure $\P_g\otimes dc$ under the mapping $(X_g,c)\mapsto X_g+c$. This measure will be understood as the "law" (it is not finite) of the field $X$ corresponding to the action \eqref{actionformal}. This measure is invariant under the shift $X\to X+a$ for any constant $a\in\R$ and is independent on the choice of $g$ conformally equivalent to $\hat g$. 

To sum up,  in what follows, we will formally understand the measure \eqref{actionformal} as the image of the product measure $\P\otimes dc$ on $H^{-1}(\R^2,\hat{g}) \times \R$ by the mapping $(X_g,c)\mapsto X_g+c$, where $dc$ is the Lebesgue measure on $\R$ and $X_g$ has the law of a GFF $X_g$ with vanishing $\lambda_g$-mean, no matter the choice of the metric $g$
conformally equivalent to $\hat g$.

\subsection{Gaussian multiplicative chaos} \label{sec:gmc}
Next we turn to the interaction term $\int e^{\gamma X}\,d\lambda_{\hat{g}}$ in eq. \eqref{actionfalse}. Since $X$ is distribution valued this is not a priori defined. As is well known it can be defined by first regularizing $X$ and then renormalizing and taking limits. This leads to the theory of Gaussian multiplicative chaos \cite{cf:Kah}. 

In what follows, we need to introduce  some cut-off approximation of the GFF $X_{g}$ for any metric $g$ conformally equivalent to the spherical metric. Natural cut-off approximations can be defined via convolution. We need that these cut-off approximations   be defined with respect to a fixed background metric: we consider Euclidean circle averages of the field because they facilitate some computations (especially Proposition \ref{circlegreen} below) but we could consider ball averages, convolutions with a smooth function or white noise decompositions of the GFF as well.

\begin{definition}{\bf (Circle average regularizations of the free field)}\label{UVcut}
We consider the field $X_{g,\epsilon}$
$$X_{g,\epsilon}(x)=\frac{1}{2\pi}\int_0^{2\pi}X_g(x+\epsilon e^{i\theta})\,d\theta.$$
\end{definition}

\begin{proposition}\label{circlegreen}
We claim (recall \eqref{theta}) 
\begin{enumerate}
\item $\lim_{\epsilon\to 0}\E[X_{\hat{g},\epsilon} (x)^2]+\ln \epsilon
+\frac{1}{2}\ln \hat{g}(x)=  \theta_{\hat{g}}+\ln 2$ uniformly on $\R^2$.
\item Let  $\psi$ be a M\"obius transform of the sphere. Denote by $(X_{\hat{g}}\circ \psi)_\epsilon$ the $\epsilon$-circle average of the field $X_{\hat{g}}\circ \psi$. Then 
$$\lim_{\epsilon\to 0}\E[(X_{\hat{g}}\circ \psi)_\epsilon(x)^2] +\frac{1}{2}\ln \hat{g}(\psi(x))+\ln|\psi'(x)|+\ln \epsilon
=\theta_{\hat{g}}+\ln 2
$$
uniformly on $\R^2$. 
\end{enumerate}
\end{proposition}

\noindent {\it Proof.}  To prove the first statement results, apply the $\epsilon$-circle average regularization to the Green function $G_{\hat{g}}$ in \eqref{covXphi} and use 
$$\int_0^{2\pi}\int_0^{2\pi}\ln\frac{1}{|e^{i\theta}-e^{i\theta'}|}\,d\theta d\theta'=0.
$$
Defining $f(x):=2m_{\hat{g}}(\ln\frac{1}{|x-\cdot|})$ and  letting  $f_\epsilon$ 
be the circle average of $f$ we then get that $\E[X_{\hat{g},\epsilon} (x)^2]+f_\epsilon(x)+ \ln \epsilon$ 
converges uniformly to $\theta_{\hat{g}}.$ Then use \eqref{second1} i.e. $f(x)=\frac{1}{2}\ln \hat{g}(x)-\ln 2$ to get the claim.

\medskip

Concerning the second statement, observe that $X_{\hat{g}}\circ \psi$ is a GFF with vanishing mean in the metric $g_\psi=|\psi'|^2\hat{g}\circ\psi$ (see Proposition \ref{GreenCC}). Therefore, the Green function of this GFF is given by \eqref{covXphi} with $g=g_\psi$. The same argument as the first item shows that 
$$
\lim_{\epsilon\to 0}(\E[X_{\hat{g}}\circ \psi_\epsilon(x)^2]+f^{\psi}_{\epsilon}(x)+\ln \epsilon)=\theta_{g_\psi}
$$
uniformly on $\R^2$ where   $f^{\psi}_{\epsilon}$ is the circle average of 
$f^{\psi}(x)=2m_{g_\psi}(\ln\frac{1}{|x-\cdot|})$. 
By \eqref{second2} 
$$f^{\psi}(x)= \frac{1}{2}\ln \hat{g}(\psi(x))+\theta_{g_\psi}+\ln|\psi'(x)|- \theta_{\hat{g}}-\ln 2$$
which yields the claim.  \qed

Define now the measure
\begin{equation}\label{Meps}
M_{\gamma,\epsilon}:=\epsilon^{\frac{\gamma^2}{2}}e^{\gamma (X_{{\hat{g}},\epsilon}+Q/2 \ln \hat{g})}  \,d\lambda
.
\end{equation}
 
\begin{proposition}\label{law}
For $\gamma\in [0,2[$,   the  following limit exists  in probability
$$
M_{\gamma}=\lim_{\epsilon\to 0}M_{\gamma,\epsilon}
=e^{\frac{\gamma^2}{2}\theta_{\hat g}+\ln 2} \lim_{\epsilon\to 0}e^{\gamma X_{\hat{g},\epsilon}-\frac{\gamma^2}{2} \E[ X_{\hat{g},\epsilon}^2] } \,d\lambda_{\hat{g}}
$$
in the sense of weak convergence of measures. This limiting measure is non trivial and is a (up to a multiplicative constant) Gaussian multiplicative chaos of the field $X_{\hat{g}}$ with respect to the measure $\lambda_{\hat{g}}$.
\end{proposition} 
 
\noindent {\it Proof.} This results  from standard tools of the general theory of Gaussian multiplicative chaos   (see \cite{review} and references therein) and Proposition \ref{circlegreen}. We also stress that all these methods were recently unified in a powerful framework in \cite{shamov}.  \qed 

\medskip
The following Proposition summarizes the behavior of this measure under M\"obius transformations:
  
\begin{proposition}\label{conformallaw}
Let $F$ be a bounded continuous function on  $H^{-1}(\R^2,\hat{g})$, $f\in \bar C({\R^2})$ and $\psi$ be  a  M\"obius transformation of the sphere. Then 
$$
(F(X_{\hat g}), \int_{\R^2} f \,dM_{\gamma})\stackrel{law}{=}(F(X_{\hat g}\circ\psi^{-1} - m_{\hat{g}_\psi}(X_{\hat g})), e^{-\gamma m_{\hat{g}_\psi}(X_{\hat g})}\int_{\R^2} f\circ\psi e^{\gamma \frac{Q}{2}\phi}dM_{\gamma})
$$
where $\hat{g}_\psi=|\psi'|^2g\circ\psi$ and $e^\phi=\hat{g}_\psi/{\hat g}$.
\end{proposition} 
\noindent {\it Proof.} 
We have
\begin{align}\nonumber
\int f \epsilon^{\frac{\gamma^2}{2}}e^{\gamma (X_{\hat{g},\epsilon}
+Q/2\ln \hat{g})}d\lambda
&=\int f\circ \psi\ \epsilon^{\frac{\gamma^2}{2}}
e^{\gamma (X_{\hat{g},\epsilon}\circ \psi+Q/2\ln \hat{g}\circ\psi )}|\psi'|^2\,d\lambda
\nonumber\\
&=\int f\circ  \psi\  (\frac{\epsilon}{|\psi'|})^{\frac{\gamma^2}{2}}
e^{\gamma (X_{\hat{g},\epsilon}\circ \psi+Q/2\ln \hat{g} )}e^{ \gamma \frac{Q}{2}\phi}\,d\lambda.
\nonumber
\end{align}
Let $\psi(z)=\frac{az+b}{cz+d}$ where $ad-bc=1$. Then $\psi'(z)=(cz+d)^{-2}$ and
$$
\phi(z)=2(\ln({1+|z|^2})-\ln({|az+b|^2+|cz+d|^2})
$$
is in $C(\overline{\R^2})$.  Let $\eta>0$. Using Proposition \ref{circlegreen} we get that
on the set $A_\eta:=B(0,\frac{1}{\eta}) \setminus B(-\frac{d}{c}, \eta)$ %
$$\lim_{\epsilon\to 0}\E[X_{\hat{g},\epsilon}(\psi(z))^2]- \E[(X_{ \hat{g}}\circ\psi)_{\frac{\epsilon}{|\psi'(z)|}}( z)^2]=0.$$
We may then use 
 the results of \cite{shamov} to conclude that the measures 
$$
(\frac{\epsilon}{|\psi'|})^{\frac{\gamma^2}{2}}
e^{\gamma (X_{\hat{g},\epsilon}\circ \psi+Q/2\ln \hat{g} )}\,d\lambda
$$
and 
$$
\epsilon^{\frac{\gamma^2}{2}} e^{\gamma (X_{\hat{g}}\circ \psi)_\epsilon+Q/2\ln \hat{g} )}\,d\lambda
$$
converge in probability to the same random measure on $A_\eta$. By Proposition
\ref{circlegreen}
\begin{equation} \label{estim2}
\E\int_{A^c_\eta}(\frac{\epsilon}{|\psi'|})^{\frac{\gamma^2}{2}}
e^{\gamma (X_{\hat{g},\epsilon}\circ \psi+Q/2\ln \hat{g} )}\lambda
\leq C\int_{A^c_\eta} ({\hat g}/\hat{g}_\psi)^{\frac{\gamma^2}{4}} \lambda_{\hat g}=
C\int_{A_\eta^c} e^{-\frac{\gamma^2}{4}\phi} \lambda_{\hat g}\to 0
\end{equation}
as $\eta\to 0$. 
By Propositions  \ref{GreenCC} and \ref{changemetric},  $X_{\hat{g}}\circ \psi$ is equal in law with  $X_{\hat{g}}-m_{\hat{g}_\psi}(X_{\hat{g}})$ yielding the claim.\qed

%

\section{Liouville Quantum Gravity on the sphere}\label{sec:LF}

In the previous Section we have given a meaning to the Gaussian part \eqref{actionformal}
of the functional integral \eqref{pathintegral} as well as for the $\mu$ term. Since the Gaussian measure we have constructed is not a finite measure one has to be careful which functionals $F$ in  \eqref{pathintegral} are integrable. Thus before giving the precise definition for  \eqref{pathintegral} we discuss this point heuristically. It turns out to lead to the first condition  in \eqref {introS}. 

As is well known for physicists the partition function i.e. the integral in \eqref{pathintegral} with $F=1$ is expected to diverge due to the integral over the  constant mode $c$ (recall that it is  distributed as the Lebesgue measure).  Let us therefore consider the toy model (sometimes called the {\it mini-superspace approximation}) where $X$ is replaced by the constant function $c$ and $DX$ by the Lebesgue measure $dc$. By  the Gauss-Bonnet theorem $\int_\Sigma R_g\lambda_g=8\pi(1-g)$ where $g$ is the genus of $\Sigma$. Therefore the
 partition function of the toy model becomes
$$\int_\R e^{-2Q(1-g)c-\mu e^{\gamma c}}\,dc.$$
This integral diverges (as $c\to -\infty$) if $g\leq 1$. Consider next $F$ to be a product of vertex operators  \eqref{vertex}. From Section \ref{sec:gmc} we know these require renormalization but proceeding within the context of the toy model they are given by $e^{\alpha_ic}$ and including them to our toy model
we end up with the integral
\begin{equation} \label{partmini}
\int_\R e^{(\sum_i \alpha_i-2Q(1-g))c-\mu e^{\gamma c}}\,dc
\end{equation}
which is finite (for $g=0$) provided the first Seiberg bound in \eqref {introS} holds. 

The divergence of the partition function has a geometric flavor. Recall that the extrema of the  Liouville action functional  \eqref{actionfalse} are  given by metrics of constant negative curvature. On the sphere and torus no such smooth metrics   exist: indeed,  it is well known that negative curvature metrics on the sphere must have conical singularities (see \cite{troyanov}). This is precisely what the vertex operators will provide via the Girsanov transform as we will now discuss.

\medskip

With these motivations we will now give the formal definition of the functional integral \eqref{pathintegral} in the presence
of the vertex operators $e^{\alpha_iX(z_i)}$. Let  $g=e^{\varphi}\hat{g}$ be a metric conformally equivalent to the spherical metric in the sense of Section \ref{sphermetric} and let $F$ be a continuous bounded functional   on $H^{-1}(\R^2,\hat{g})$. We define
\begin{align}\label{eq:metricg}
\Pi_{\gamma,\mu}^{(z_i\alpha_i)_i}& (g,F;\epsilon)\\
 :=& e^{\frac{1}{96\pi}\int_{\R^2}|\partial^{\hat{g}} \varphi|^2+2R_{\hat{g}}\varphi \,d\lambda_{\hat{g}}}\int_\R\E\Big[F( c+X_{g}+Q/2\ln g)\prod_i \epsilon^{\frac{\alpha_i^2}{2}}e^{\alpha_i (c+X_{g,\epsilon}+Q/2\ln g)(z_i)}\nonumber\\
&\exp\Big( - \frac{Q}{4\pi}\int_{\R^2}R_{g} (c+X_{g })  \,d\lambda_{g} - \mu 
\epsilon^{\frac{\gamma^2}{2}}\int_{\R^2}e^{\gamma (c+X_{g,\epsilon}+ Q/2\ln g) }\,d\lambda \Big) \Big]\,dc.\nonumber %
\end{align}
and we want to inquire when  the limit 
$$\lim_{\epsilon\to 0}\Pi_{\gamma,\mu}^{(z_i\alpha_i)_i} (g,F;\epsilon)=:\Pi_{\gamma,\mu}^{(z_i\alpha_i)_i} (g,F)
$$ 
exists.

\begin{remark} We include 
the additional factor $e^{\frac{1}{96\pi}\int_{\R^2}|\partial^{\hat{g}} \varphi|^2+2R_{\hat{g}}\varphi \,d\lambda_{\hat{g}}}$ to conform to the physics conventions. Indeed  the formal expression \eqref{pathintegral} differs from \eqref{eq:metricg} in that in the latter we use a normalized expectation for the Free Field. Thus to get \eqref{pathintegral} we would need to multiply by the Free Field partition function $z(g)$. The latter is not uniquely defined but its variation with metric is:
$$
z(e^\varphi \hat g)=e^{\frac{1}{96\pi}\int_{\R^2}|\partial^{\hat{g}} \varphi|^2+2R_{\hat{g}}\varphi \,d\lambda_{\hat{g}}}z(\hat g)
$$
 see \cite{dubedat,gaw}. This additional factor makes the Weyl anomaly formula conform with the 
 standard one in Conformal Field Theory. We note also that
 the translation by $Q/2\ln g$ in the argument of $F$ is necessary for conformal invariance (Section \ref{Conformal covariance}).
\end{remark}

We start by considering the round metric, $g=\hat g$.  We first handle the curvature term. Since $R_{\hat g}=2$ and  $X_{\hat{g}}$ has vanishing $\lambda_{\hat{g}}$-mean we obtain 
\begin{align}\label{eq:defPiGB}
\Pi_{\gamma,\mu}^{(z_i\alpha_i)_i}& (\hat{g},F;\epsilon)\\
 =&  \int_\R e^{-2Qc}\,\E\Big[F( c+X_{\hat{g} } +Q/2\ln\hat{g})\prod_i \epsilon^{\frac{\alpha_i^2}{2}}e^{\alpha_i (c+X_{\hat{g},\epsilon}+Q/2\ln\hat{g})(z_i)}\nonumber\\
&\exp\Big( - \mu \epsilon^{\frac{\gamma^2}{2}}\int_{\R^2}e^{\gamma (c+X_{\hat{g},\epsilon}+ Q/2\ln\hat{g} )}\,d\lambda \Big) \Big]\,dc
.\nonumber %
\end{align}

Now we handle the insertions operators $e^{\alpha_i  X_{\hat{g},\epsilon}(z_i)}$. In view of  Proposition \ref{circlegreen}, we can write (with the Landau notation)
\begin{equation}\label{def:insert}
\epsilon^{\frac{\alpha_i^2}{2}}e^{\alpha_i X_{\hat{g},\epsilon}(z_i)}=e^{\frac{\alpha_i^2}{2}(\theta_{\hat{g}}+\ln 2)}\hat{g}(z_i)^{-\frac{\alpha_i^2}{4}}e^{\alpha_i X_{\hat{g},\epsilon}(z_i)- \frac{\alpha_i^2}{2}\E[X_{\hat{g},\epsilon}(z_i)^2]}(1+o(1)).
\end{equation}
Note that the $o(1)$ term is deterministic as it just comes from the normalization of variances. Then, by applying the Girsanov transform and setting
\begin{equation}\label{def:Hepsilon}
H_{\hat{g},\epsilon}(x)=\sum_i\alpha_i\int_0^{2\pi}G_{\hat{g}}(z_i+\epsilon e^{i\theta},x)\frac{d\theta}{2\pi},
\end{equation}
we obtain  
\begin{align}\label{eq:insert}
\Pi_{\gamma,\mu}^{(z_i\alpha_i)_i}& (\hat{g},F;\epsilon)
 =e^{C_{\epsilon}({\bf z})}\Big(\prod_i\hat{g}(z_i)^{-\frac{\alpha_i^2}{4}+\frac{Q}{2}\alpha_i}\Big) \\&\int_\R e^{\big(\sum_i\alpha_i-2Q\big)c}\,\E\Big[F( c
 +X_{\hat{g}} +H_{\hat{g},\epsilon}+Q/2\ln\hat{g}) (1+o(1))\nonumber \\&\times  
  \exp\Big( - \mu e^{\gamma c}\epsilon^{\frac{\gamma^2}{2}}\int_{\R^2}e^{\gamma( X_{\hat{g},\epsilon}+H_{\hat{g},\epsilon}+ Q/2\ln\hat{g}) }\,d\lambda \Big) \Big]\,dc,\nonumber %
\end{align}
with
\begin{align}\label{def:chatg}
\lim_{\epsilon\to 0}C_{\epsilon}({\bf z})=&\frac{1}{2}  \sum_{i\not =j}\alpha_i\alpha_jG_{\hat{g}}(z_i,z_j)
+\frac{\theta_{\hat{g}}+\ln 2}{2}\sum_i\alpha_i^2:=C({\bf z}).
\end{align}
In the next subsection we study under what conditions the limit of \eqref{eq:insert} exists.

\subsection{Seiberg bounds and KPZ scaling laws}\label{sub:seiberg}
Since $H_{\hat{g},\epsilon}$ converges in $H^{-1}(\R^2,\hat{g})$ to
\begin{equation}\label{def:H}
H_{\hat{g}}(x)=\sum_i\alpha_iG_{\hat{g}}(z_i,x)
\end{equation}
it suffices to study the convergence of the partition function $\Pi_{\gamma,\mu}^{(z_i\alpha_i)_i} (\hat{g},1;\epsilon)$. We show that 
 a necessary and sufficient condition for the Liouville partition function to have a non trivial 
 limit is the validity of the { Seiberg bounds} given in eq. \eqref {introS}.

It is in fact easy to see that the first Seiberg bound is a necessary condition even for the existence of the regularized theory. Indeed, let 
\begin{equation}\label{def:Zeps}
Z_\epsilon:=\epsilon^{\frac{\gamma^2}{2}}\int_{\R^2}e^{\gamma( X_{\hat{g},\epsilon}+H_{\hat{g},\epsilon}+ Q/2\ln\hat{g}) }\,d\lambda.
\end{equation}
Note that $|H_{\hat{g},\epsilon}(z)|\leq C_\epsilon$ since $G(z_i,z)$ tends 
to constant as $|z|\to\infty$. Hence from Proposition \ref{law} we infer $\E [Z_\epsilon]<\infty$
and thus $Z_\epsilon<\infty$  $\P$-almost surely. Hence we can find $A>0$ such that $\P(Z_\epsilon\leq A)>0$ and then
\begin{align*}
\Pi_{\gamma,\mu}^{(z_i\alpha_i)_i}  (\hat{g},1,\epsilon)\geq \Big(\prod_i\hat{g}(z_i)^{-\frac{\alpha_i^2}{4}+\frac{Q}{2}\alpha_i}\Big) e^{C_{\epsilon}({\bf z})}\int_{-\infty}^0e^{\big(\sum_i\alpha_i-2Q\big)c}e^{-\mu e^{\gamma c}A}\P(Z_\epsilon\leq A)\,dc=+\infty
\end{align*}
if the first condition in  \eqref {introS}  fails to hold.  We will see shortly that the second condition 
 $\alpha_i<Q$ is needed to ensure that the integral in \eqref{def:Zeps}
 does not blow up in the neighborhood of the places of insertions $(z_i)_i$ as 
 $\epsilon \to 0$.

Finally, we mention that the bounds  \eqref {introS}  show that the number of vertex operator insertions must be at least  $3$ in order to have well defined correlation functions of the Liouville theory on the sphere. This conforms with the fact  (see \cite{troyanov}) that on the sphere one must  insert at least three conical singularities in order to construct a metric with negative curvature. We claim 
\begin{theorem}{\bf (Convergence of the partition function)}\label{th:seiberg} 
Let $\sum_i\alpha_i>2Q$.   Then the limit 
$$\lim_{\epsilon\to 0}\Pi_{\gamma,\mu}^{(z_i\alpha_i)_i}  (\hat{g},1;\epsilon):=\Pi_{\gamma,\mu}^{(z_i\alpha_i)_i}  (\hat{g},1)$$
     exists. The limit
is nonzero if $\alpha_i< Q$ for all $i$ whereas it vanishes identically if  $\alpha_i\geq Q$ for some $i$. 
\end{theorem}
  \noindent  {\it Proof.} Eq.
 \eqref{eq:insert} gives for $F=1$
\begin{align}\nonumber
\Pi_{\gamma,\mu}^{(z_i\alpha_i)_i}  (\hat{g},1,\epsilon)  =\prod_i\hat{g}(z_i)^{-\frac{\alpha_i^2}{4}+\frac{Q}{2}\alpha_i}e^{C({\bf z})}(1+o(1))
  \E\Big[\int_\R e^{c(\sum_i\alpha_i-2Q)}\exp\big(    - \mu e^{\gamma c}Z_\epsilon\big) \,dc\Big] .
  \end{align}
 As remarked above,  $Z_\epsilon>0$ almost surely.  
By making the change of variables $u=\mu e^{\gamma c}Z_\epsilon$ in \eqref{eq:insert}, we compute 
\begin{align}\label{int1}
 \E\Big[\int_\R e^{c(\sum_i\alpha_i-2Q)}\exp\big(    - \mu e^{\gamma c}Z_\epsilon\big) \,dc\Big] 
= \frac{\mu^{\frac{\sum_i\alpha_i-2Q}{\gamma}}}{\gamma}\Gamma\Big(\gamma^{-1}(\sum_i\alpha_i-2Q)\Big)\E\Big[\frac{1}{Z_\epsilon^{\frac{\sum_i\alpha_i-2Q}{\gamma}}}\Big]
\end{align} 
where $\Gamma$ is the standard $\Gamma$ function.  The claim follows from the following Lemma.
\qed
\begin{lemma}\label{exist} Let $s<0$.
 If $\alpha_i<Q$ for all $i$ then 
$$
\lim_{\epsilon\to 0}\E [Z_\epsilon^s]=\E [Z_0^s]
$$
where 
\begin{equation}\label{Z0def}
Z_0=\int_{\R^2} e^{\gamma H_{\hat g}(x)}
M_{\gamma}(dx)
\end{equation}
 and the limit is nontrivial: $0<\E Z_0^s<\infty$. 
  
  If $\alpha_i\geq Q$ for some $i\in\{1,\dots, p\}$ then 
$$ \lim_{\epsilon\to 0}\E Z_\epsilon^s= 0.$$
\end{lemma}
\medskip

As a corollary of the relation \eqref{int1}, we obtain a rigorous derivation of the KPZ 
scaling laws (see \cite{cf:KPZ,nakayama,witten} for physics references)
\begin{theorem}{\bf (KPZ scaling laws)}\label{KPZscaling}
We have the following exact scaling relation for the Liouville partition function  with insertions $(z_i,\alpha_i)_i$
$$\Pi_{\gamma,\mu}^{(z_i\alpha_i)_i}  (\hat{g},1) =\mu^{\frac{2Q-\sum_i\alpha_i}{\gamma}}\Pi_{\gamma,1}^{(z_i\alpha_i)_i}  (\hat{g},1)  
$$ 
where 
$$
\Pi_{\gamma,1}^{(z_i\alpha_i)_i}  (\hat{g},1)  =e^{C({\bf z})}\Big(\prod_i\hat{g}(z_i)^{\Delta_{\alpha_i}}\Big) \gamma^{-1}\Gamma\Big(\gamma^{-1}(\sum_i\alpha_i-2Q)\Big)\E\Big[\frac{1}{Z_0^{\frac{\sum_i\alpha_i-2Q}{\gamma}}}\Big]
$$ 
and we defined
\begin{equation}\label{weight}
\Delta_{\alpha}=\frac{\alpha}{2}\big(Q-\frac{\alpha}{2}\big)
\end{equation}
and $C({\bf z})$ is defined by \eqref{def:chatg}.  Moreover
\begin{align}\label{eq:insert1}
\Pi_{\gamma,\mu}^{(z_i\alpha_i)_i}& (\hat{g},F)
 =e^{C({\bf z})}\prod_i\hat{g}(z_i)^{\Delta_{\alpha_i}}
 \\&\int_\R e^{\big(\sum_i\alpha_i-2Q\big)c}\,\E\Big[F( c
 +X_{\hat{g}} +H_{\hat{g}}+Q/2\ln\hat{g}) \exp\Big( - \mu e^{\gamma c}Z_0 \Big) \Big]\,dc.\nonumber %
\end{align}

\end{theorem}

\medskip 

\noindent {\it Proof of Lemma \ref{exist}.} Note first that $\E Z_\epsilon^s<\infty$ for
all $\epsilon\geq 0$. Indeed, 
recalling \eqref{Meps} 
$$
Z_\epsilon=
\int_{\R^2} e^{\gamma H_{\hat{g},\epsilon}(z)}M_{\gamma,\epsilon}(dx).
$$
Take any non empty ball $B$ that contains  no $z_i$. Then
$$\E [Z_\epsilon^s]\leq A^s\E  [M_{\gamma,\epsilon}(B)^s]
$$
where  $A=C\min_{z\in B}\frac{4e^{ \gamma H_{\hat{g}}(z)}}{(1+|z|^2)^2}$. It is a standard fact in Gaussian multiplicative chaos theory (see \cite[Th 2.12]{review} again) that the random variable $M_{\gamma,\epsilon}(B)$ possesses negative moments of all orders for $\gamma\in [0,2[$.

Let now $\alpha_i<Q$ for all $i$.
 Let us consider the set
$A_r=\cup_i B(z_i,r)$ and write
$$Z_\epsilon=\int_{A_r} e^{\gamma H_{\hat{g},\epsilon}(z)}M_{\gamma,\epsilon}(dx)+
\int_{A_r^c} e^{\gamma H_{\hat{g},\epsilon}(z)}M_{\gamma,\epsilon}(dx):=Z_{r,\epsilon}+Z_{r,\epsilon}^c.
$$ 
Since $H_{\hat{g},\epsilon}$ converge uniformly on $A_r^c$  to a continuous limit
the limit
\begin{equation}\label{nonsingular}
\lim_{\epsilon\to 0}Z_{r,\epsilon}^c=\int_{A_r^c} e^{\gamma H_{\hat{g,\epsilon}}(z)}M_{\gamma}(dx):=Z_{r,0}^c
\end{equation}
exists in probability by Proposition \ref{law}. 

We study next the $r$-dependence of $Z_{r,\epsilon}$. Without loss of generality, we may take $\epsilon=2^{-n}$
and $r=2^{-m}$ with $n>m$ and $A_r=B(0,r)$. Then, dividing  $B(0,r)$ to
dyadic annuli $2^{-k-1}\leq |z|\leq 2^{-k}$ and noting that  $e^{\gamma H_{\hat{g},\epsilon}(z)}\leq C2^{\gamma\alpha k}$ on such annulus we get 
\begin{equation}\label{ZREPS}
Z_{r,\epsilon}=\int_{B(0,r)} e^{\gamma H_{\hat{g},\epsilon}(z)}M_{\gamma,\epsilon}(dx)\leq 
C\sum_{k=m}^{n} 2^{\gamma\alpha k}M_{\gamma,\epsilon}(B_k)
\end{equation}
where $B_k=B(0,2^{-k})$. 

The distribution of $M_{\gamma,\epsilon}(B_k)$ is easiest to study using  the white noise cutoff $(\tilde X_\epsilon)_\epsilon$ of $X_{\hat{g}}$. More precisely, the family $(\tilde X_\epsilon)_\epsilon$ is a family of Gaussian processes defined as follows. Consider the heat kernel $(p_t(\cdot,\cdot))_{t\geq 0}$ of the Laplacian $\triangle_{\hat{g}}$ on $\R^2$. Let $W$ be a white noise distributed  on $\R_+\times\R^2$ with intensity $dt\otimes \lambda_{\hat{g}}(dy)$. Then
$$\tilde X_{\epsilon} (x)=\frac{1}{\sqrt{2\pi}}\int_{\epsilon^2}^{\infty}\big(p_{t/2}(x,y)-\frac{1}{\lambda_{\hat{g}}(\R^2)}\big)W(dt,dy).$$
The correlation structure of the family $(\tilde X_{\epsilon} )_{\epsilon>0}$ is given by
\begin{equation}\label{def:covWN}
\E[\tilde X_{\epsilon} (x)\tilde X_{\epsilon'}  (x')]=\frac{1}{2\pi}\int_{(\epsilon\wedge \epsilon')^2}^{\infty}\big(p_{t} (x,x')-\frac{1}{\lambda_{\hat{g}}(\R^2)}\big) \,dt.
\end{equation}
For $\epsilon>0$, we define the random measure
$$\tilde M_{\gamma,\epsilon}:= e^{\gamma \tilde X_{\epsilon}-\frac{\gamma^2}{2}\E[(\tilde X_{\epsilon}(x))^2]}  \,d\lambda_{\hat{g}}$$ and 
$\tilde M_{\gamma}:=\lim_{\epsilon\to 0}\tilde M_{\gamma,\epsilon}$, which has the same law as $M_\gamma$ (see \cite[Thm 3.7]{review}).
The  covariance  of the field  $X_{\hat g,\epsilon}$ is comparable to the one of 
$\tilde X_{\epsilon}$. Indeed, uniformly in $\epsilon$,
$$
\E[\tilde X_{\epsilon} (x)\tilde X_{\epsilon}  (y)]\leq C+\E[X_{\epsilon} (x) X_{\epsilon}  (y)]
$$
and so by  Kahane's convexity inequality (see \cite{cf:Kah}) we get, for $q\in (0,1)$
$$
\E[M_{\gamma,\epsilon}(B_k)^q]\leq C\E[\tilde M_{\gamma,\epsilon}(B_k)^q].
$$
We have the relation
\begin{equation}\label{moments}
\sup_{\epsilon}\E[\tilde M_{\gamma,\epsilon}(B_k)^q]\leq  C_q2^{-k\xi(q)}
\end{equation}
for all $q<\frac{4}{\gamma^2}$  where $\xi(q)=(2+\frac{\gamma^2}{2})q-\frac{\gamma^2}{2}q^2$. Indeed, the family $(\tilde M_{\gamma,\epsilon}(B_k))_{\epsilon}$ is a martingale so that, by Jensen, it suffices to prove that the limit $\tilde M_{\gamma }$ satisfies such a bound. This latter fact is standard, see \cite[Th 2.14]{review} for instance.

Therefore by Tchebychev
$$
\P(Z_{r,\epsilon}>R) 
\leq C_{q,\delta}R^{-q}\sum_{k=m}^n 2^{-k\xi(q)}2^{(\gamma\alpha+\delta)q k}\leq
C_{q,\delta}R^{-q}2^{-m(\xi(q)-q(\gamma\alpha+\delta))}
$$
provided $(\gamma\alpha+\delta)q<\xi(q)$.  This holds for $q$ and $\delta$ small enough since $\alpha<Q$ i.e.
$\gamma\alpha< 2+\frac{\gamma^2}{2}$.  Hence, for some $\alpha,\beta>0$
\begin{equation}\nonumber
\P(Z_{r,\epsilon}>r^\alpha) \leq Cr^\beta \ \ \ \forall\epsilon\geq 0
\end{equation}
where we noted that the same argument covers also the $\epsilon=0$ case. 

Let $\chi_r=1_{Z_{r,\epsilon}>r^\alpha}$. 
We get  by Schwartz
$$
|\E[((Z_{r,\epsilon}+Z_{r,\epsilon}^c)^s-(Z_{r,\epsilon}^c)^s)\chi_r]|\leq 2(\E\chi_r\E (Z_{r,\epsilon}^c)^{2s})^{1/2}\leq Cr^{\beta/2}(\E (Z_{r,\epsilon}^c)^{2s})^{1/2}
$$
and  using $|(a+b)^s-b^s|\leq Cab^{s-1}$
$$
|\E((Z_{r,\epsilon}+Z_{r,\epsilon}^c)^s-(Z_{r,\epsilon}^c)^s)(1-\chi_r)|\leq Cr^\alpha
\E(Z_{r,\epsilon}^c)^{s-1}.
$$
Since $\E (Z_{r,\epsilon}^c)^{s}\leq \E (Z_{1,\epsilon}^c)^{s}$ and the latter stays bounded
as $\epsilon\to 0$ we conclude 
$$
|\E[(Z_{\epsilon})^s-(Z_{r,\epsilon}^c)^s)]|\leq C (r^\alpha+r^\beta)
$$
for all $\epsilon\leq r$. In particular, for $\epsilon=0$ this gives   
\begin{equation}\label{PZ_r}
\lim_{r\to 0}\E[(Z_{r,0}^c)^s]=\E [Z_0^s].
 \end{equation}
Since  $\E [(Z_{r,\epsilon}^c)^s]<\infty$ for
all $\epsilon\geq 0$  and by \eqref{nonsingular}  $Z_{r,\epsilon}^c$ converges in probability to $Z_{r,0}^c$ as $
\epsilon\to 0$ we have $\lim_{\epsilon\to 0}\E [(Z_{r,\epsilon}^c)^s]=\E [(Z_{r,0}^c)^s]$. From
\eqref{PZ_r} we then conclude our claim
$\lim_{\epsilon\to 0}\E[(Z_{\epsilon})^s]=\E[(Z_{0})^s]$. 

For later purpose let us remark that from \eqref{moments} we get
$$
M_\gamma(B_k)\leq C_\delta(\omega)2^{-k(2+\frac{\gamma^2}{2}-\delta)}
$$
where $C_\delta(\omega)<\infty$ almost surely. This easily leads to
\begin{equation}\label{anotherlimit}
\sup_{\epsilon>0}\int_{B_r} e^{\gamma H_{\hat{g},\epsilon}(z)}M_{\gamma}(dx)\to 0
\end{equation}
in probability as $r\to 0$.

\medskip
Let us now prove the second part of the lemma. Without loss of generality, we may assume that $\alpha_1\geq Q$ and $z_1=0$. It suffices to prove for
 the $Z_{1,\epsilon}$ defined in \eqref{ZREPS} that
\begin{equation}\label{ZREPS1}
\lim_{\epsilon\to 0}\E[ Z_{1,\epsilon}^s]= 0
\end{equation}
By Kahane convexity \cite{cf:Kah} (or \cite[Thm 2.1]{review}) we get
$$
\E[ Z_{1,\epsilon}^s]\leq C\E[ \tilde Z_{1,\epsilon}^s].
$$
Next, we  bound 
\begin{equation}\label{z1ep}
\tilde Z_{1,\epsilon}\geq c\sum_{k=1}^n 2^{\alpha\gamma  k}\tilde M_{\gamma,\epsilon}(A_k)
\geq c\max_{k\leq n}2^{(2+\gamma^2/2)k}\tilde M_{\gamma,\epsilon}(A_k)
\end{equation}
where $A_k$ is the annulus with radi $2^{-k}$ and $2^{-k+1}$ and we recall
that $\epsilon=2^{-n}$ and $\alpha\gamma\geq 2+\gamma^2/2$ . We may then decompose, for $r=2^{-k}$ (and $\epsilon<r$),
\begin{equation}\label{decomp1}
\tilde M_{\gamma,\epsilon}(dz)=e^{\gamma \tilde X_r(z)-\frac{\gamma^2}{2}\E[\tilde X_r(z)^2]}r^2\widehat{M}_{\gamma,\epsilon, r}(dz/r)
\end{equation} 
where the measure $\widehat{M}_{\gamma,\epsilon, r}$ is independent of the sigma-field $\{\tilde X_u(x);u\geq r,x\in\R^2\}$ and has the law
$$\widehat{M}_{\gamma,\epsilon, r}(dz)=e^{\gamma (\tilde X_{\epsilon}-\tilde X_r)(rz)-\frac{\gamma^2}{2}\E[(\tilde X_{\epsilon}-\tilde X_r)(rz)^2]}\,dz.$$
We can rewrite \eqref{decomp1} as 
\begin{equation}\label{decomp2}
\tilde M_{\gamma,\epsilon}(dz)=e^{\gamma \tilde X_r(0)-\frac{\gamma^2}{2}\E[(\tilde X_r(0))^2]}e^{\gamma (\tilde X_r(z)-\tilde X_r(0))-\frac{\gamma^2}{2}(\E[(\tilde X_r(z))^2]-\E[(\tilde X_r(0))^2])}r^2\widehat{M}_{\gamma,\epsilon, r}(dz/r)
\end{equation} 
to get
\begin{equation}\label{decomp3}
\tilde M_{\gamma,\epsilon}(A_k)\geq r^2 e^{\gamma \tilde X_r(0)-\frac{\gamma^2}{2}\E[(\tilde X_r(0))^2]}e^{\min_{z\in B(0,1)} Y_r(z)}\widehat{M}_{\gamma,\epsilon, r}(A_1)
\end{equation} 
with $Y_r(z)=\gamma(\tilde X_r(rz)-\tilde X_r(0))-\frac{\gamma^2}{2}(\E[(\tilde X_r(rz))^2]-\E[(\tilde X_r(0))^2])$. Now we want to determine the behavior of all the terms involved in the above right-hand side.

By using in turn Doob's inequality and then Kahane convexity \cite{cf:Kah} (or \cite[Thm 2.1]{review}), we get 
\begin{equation}\label{momentq}
\E [\sup_{\epsilon<r}\widehat{M}_{\gamma,\epsilon, r}(A_1)^{-q}]\leq c_q\E [\widehat{M}_{\gamma,0, r}(A_1)^{-q}]\leq \E[M_{\gamma}(A_1)^{-q}]\leq C_q.
\end{equation}
uniformly in $r\leq 1$. Hence, for all $a>0$
\begin{equation}\label{momentq1}
\P (\sup_{\epsilon<r} \widehat{M}_{\gamma,\epsilon, r}(A_1)\leq n^{-1})\leq  C_a n^{-a}.
\end{equation}

Next, we estimate the $\min$ in \eqref{decomp3}. The key point is to observe that the Gaussian process $Y_r$ does not fluctuate too much in such a way that its minimum possesses a   Gaussian left tail distribution. To prove this, we write $Y_r(z)=\E [Y_r(z)]+Y'_r(z)$ and we note that using the covariance structure of $(\tilde X_r)_r$
we get for all $z\in B(0,1)$
$$
|\E Y_r(z)|=\frac{\gamma^2}{2}|\E[(\tilde X_r(rz))^2]-\E[(\tilde X_r(0))^2]|\leq C
$$
and for all $ z,z'\in B(0,1)$, 
$$   \E[(Y'_r(z) -Y'_r(z') )^2]\leq C|z-z'|,
$$
uniformly in $r\leq 1$. Using for example \cite[Thm. 7.1, Eq. (7.4)]{Ledoux}, one can then deduce  
$$
\forall x\geq 1,\quad \sup_r\P (\min_{z\in B(0,1)}\gamma Y_r(z)\leq -x)\leq Ce^{-cx^2}
$$ 
for some constants $C,c>0$. Hence, for all $a>0$
\begin{equation}\label{momentq2}
\P ( e^{\min_{z\in B(0,1)} Y_r(z)}\leq n^{-1})\leq  C_a n^{-a}.
\end{equation}
Combining \eqref{decomp3}, \eqref{momentq1} and \eqref{momentq2} with \eqref{z1ep} we conclude
$$
\P(\tilde Z_{1,\epsilon}<n)\leq \P(\max_{k\leq n}e^{\gamma X_{2^{-k}}(0)}\leq n^3 )+Cn^{-a}.
$$
Since the law of the path $t\mapsto \tilde X_t(0)$   is that of  Brownian motion at time $-\ln t$ the first term on the RHS  
tends to zero as $n\to\infty$ and \eqref{ZREPS1} follows.\qed

\subsection{Conformal covariance, KPZ formula and Liouville field}\label{Conformal covariance}

In what follows, we assume that the bounds  \eqref {introS} hold and we will study how the $n$-point correlation functions  $\Pi_{\gamma,\mu}^{(z_i\alpha_i)_i} (\hat{g},F)$ transform under conformal reparametrization of the sphere. The KPZ formula describes precisely the rule for these transformations.   We claim (recall \eqref{weight})

\begin{theorem}{\bf (Field theoretic KPZ formula)}\label{KPZ}
Let  $\psi$ be a M\"obius transform  of the sphere. Then
$$\Pi_{\gamma,\mu}^{(\psi(z_i),\alpha_i)_i} (\hat{g},1)=\prod_i|\psi'(z_i)|^{-2\Delta_{\alpha_i}}\Pi_{\gamma,\mu}^{(z_i,\alpha_i)_i} (\hat{g},1).$$
\end{theorem}

Let us now define the law of the Liouville field on the sphere.
\begin{definition}{\bf (Liouville field)}\label{def:LF}
We define a probability law $\P^{  \gamma,\mu}_{(z_i,\alpha_i)_i,\hat{g}}$ on $H^{-1}(\R^2,\hat{g})$(with expectation $\E^{  \gamma,\mu}_{(z_i,\alpha_i)_i,\hat{g}}$) by
\begin{align*}
\E^{  \gamma,\mu}_{(z_i,\alpha_i)_i,\hat{g}}[F(\phi)]=\frac{\Pi_{\gamma,\mu}^{(z_i\alpha_i)_i} (\hat{g},F)}{\Pi_{\gamma,\mu}^{(z_i\alpha_i)_i} (\hat{g},1)},
\end{align*}
for all bounded continuous functional on $H^{-1}(\R^2,\hat{g})$.  
\end{definition}

We have the following result about the behaviour of the Liouville field under the M\"obius transforms of the sphere
\begin{theorem}\label{coroCR}
Let $\psi$ be a M\"obius transform of the sphere. The law of the Liouville field $\phi$ under $\P^{  \gamma,\mu}_{( z_i,\alpha_i)_i,\hat{g}}$ is the same as that of $\phi\circ \psi+Q\ln|\psi'|$ under $\P^{  \gamma,\mu}_{(\psi(z_i),\alpha_i)_i,\hat{g}}$.
 \end{theorem}

\noindent {\it Proof of Theorems \ref{KPZ} and   \ref{coroCR}.} 
We start from the relation  \eqref{eq:insert1}. Let 
$$H_{\hat g}^\psi(z)=\sum_i\alpha_iG_{\hat{g}}(\psi(z_i),z).
$$ We apply Proposition
\ref{conformallaw} to $f= e^{\gamma H^\psi_{\hat g,\epsilon}}$. By \eqref{anotherlimit} we can take the limit  $\epsilon\to0$
 to get
\begin{align}\nonumber
\Pi_{\gamma,\mu}^{(\psi(z_i),\alpha_i)_i}& (\hat{g},F)
 =e^{C({\bf \psi(z)})}\prod_i\hat{g}(\psi(z_i))^{\Delta_{\alpha_i}}
 \int_\R e^{sc}
\E\Big[
F( c +X_{\hat{g}}\circ\psi^{-1} -m_{\hat{g}_\psi}(X_{\hat g})+H^\psi_{\hat{g}}+Q/2\ln\hat{g})\\&
  \exp\big(- \mu e^{\gamma (c-m_{\hat{g}_\psi}(X_{\hat g}))} 
 \int  e^{\gamma( H^\psi_{\hat{g}}\circ\psi+
  \frac{Q}{2}\phi)}dM_{\gamma}\big)
]\,dc.\nonumber %
\end{align}
where we denoted $s=\sum_i\alpha_i-2Q$.
Next, use  the shift invariance of the Lebesgue measure (we make the change of variables $c=c'+m_{\hat{g}_\psi}(X_{\hat{g}})$)  
to get  
\begin{align}\label{eq:insert2}
\Pi_{\gamma,\mu}^{(\psi(z_i)\alpha_i)_i}  (\hat{g},F)
 =&e^{C({\bf \psi(z)})}\prod_i\hat{g}(\psi(z_i))^{\Delta_{\alpha_i}}
\int_\R e^{sc}
\E \Big[e^{sm_{\hat{g}_\psi}(X_{\hat g})}
 F( c+X_{\hat{g}}\circ\psi^{-1} +H_{\hat{g},\psi}+Q/2\ln\hat{g})\\&
  \exp\big(- \mu e^{\gamma c}\int  e^{\gamma( H_{\hat{g},\psi}\circ\psi+ \frac{Q}{2}\phi)}dM_{\gamma}\big)
\Big]\,dc.\nonumber %
\end{align}
Now we apply the Girsanov transform to the term $e^{sm_{\hat{g}_\psi}(X_{\hat g})}$ where $m_{\hat{g}_\psi}(X_{\hat g})=\frac{1}{4\pi}\int X_{\hat g}e^\phi d\lambda_{\hat g}$ and $e^\phi =\frac{|\psi'|^2\hat{g}\circ\psi}{\hat{g}}$. This has the effect of shifting the law of the field $X_{\hat{g}}$, which becomes 
$$X_{\hat{g}}+\frac{s}{4\pi}G_{\hat{g}}e^\phi.$$
The variance of this Girsanov transform is 
${s^2}D_\psi $ where 
\begin{equation}\label{Cpsi'}
D_\psi= \frac{1}{4\pi}m_{\hat g}(e^\phi G_{\hat{g}}e^\phi)=\frac{1}{(4\pi)^2}
\int_{\R^2}\int_{\R^2}  G_{\hat{g}}(z,z') \lambda_{{g_\psi}}(dz)\lambda_{{g_\psi}}(dz'), 
\end{equation}
i.e. the whole partition function will be multiplied by $e^{\frac{s^2}{2} D_\psi   }$. 

Plugging in the shifted field to \eqref{eq:insert2} we need to compute $
H_{\hat{g},\psi}\circ\psi+\frac{s}{4\pi}G_{\hat{g}}e^\phi$. First, using \eqref{Grule} for
$ (H^\psi_{\hat{g}}\circ\psi)(z)=\sum_i\alpha_iG_{\hat g}(\psi(z),\psi(z_i))$ we get
$$
H^\psi_{\hat{g}}\circ\psi=H_{\hat{g}}-\frac{\sum\alpha_i}{4}\phi(z)-\frac{1}{4}\sum_i\alpha_i\phi(z_i).
$$
Next, to compute $G_{\hat{g}}e^\phi$ note that both metrics $\hat g$ and $\hat g_\psi=e^\phi\hat g$
have Ricci curvature $2$. Hence from \eqref{curvature} we infer $e^\phi=1-\frac{1}{2}\Delta_{\hat g}\phi$ and thus 
\begin{equation}\label{Gonephi}
\frac{1}{4\pi}G_{\hat{g}}e^\phi=\frac{1}{4}(\phi-m_{\hat g}(\phi)).
\end{equation}
Combining we get
$$
H_{\hat{g},\psi}\circ\psi+\frac{s}{4\pi}G_{\hat{g}}e^\phi=H_{\hat{g}}-\frac{Q}{2}\phi(z)-\frac{1}{4}\sum_i\alpha_i\phi(z_i)-\frac{s}{4}m_{\hat g}(\phi).
$$
Thus  \eqref{eq:insert2} becomes 
\begin{align}\nonumber 
\Pi_{\gamma,\mu}^{(\psi(z_i)\alpha_i)_i}& (\hat{g},F)
 =e^{C({\bf\psi(z)})}\Big(\prod_i\hat{g}(\psi(z_i))^{\Delta_{\alpha_i}}
\Big) \int_\R e^{sc}
\E 
\Big[F\big( c'
 +(X_{\hat{g}}+H_{\hat{g}}+Q/2(\ln\hat{g}-\ln|\psi'|^2))\circ\psi^{-1}\big)\\&
  \exp\big(- \mu e^{\gamma c'}\int  e^{\gamma H_{\hat{g}}}dM_{\gamma}\big)
]\,dc\  e^{\frac{s^2}{2} D_\psi   }.\nonumber %
\end{align}
where
$$
c'=c-\frac{s}{4}m_{\hat g}(\phi)
-\frac{1}{4}\sum_i\alpha_i\phi(z_i).
$$
By a shift in the $c$-integral we get
\begin{align}\nonumber 
\Pi_{\gamma,\mu}^{(\psi(z_i)\alpha_i)_i}& (\hat{g},F)
 =e^{C({\bf \psi(z)})}\prod_i\hat{g}(\psi(z_i))^{\Delta_{\alpha_i}}
 \int_\R e^{sc}
\E 
\Big[F\big( c
 +(X_{\hat{g}}+H_{\hat{g}}+Q/2(\ln\hat{g}-\ln|\psi'|^2))\circ\psi^{-1}\\&
  \exp\big(- \mu e^{\gamma c}\int  e^{\gamma H_{\hat{g}}}dM_{\gamma}\big)
]\,dc\  e^{\frac{s}{4}\sum_i\alpha_i\phi(z_i)}e^{\frac{s^2}{2} (D_\psi +\frac{1}{2}m_{\hat g}(\phi) ) }%
\end{align}
Combining  \eqref{def:chatg} with \eqref{Grule} we have 
$$
C({\bf \psi(z)})=C({\bf z})-\frac{1}{8}\sum_{i\not =j}\alpha_i\alpha_j(\phi(z_i)+\phi(z_j))=C({\bf z})
-\frac{\sum_{i}\alpha_i}{4}\sum_j\alpha_j\phi(z_j)+\frac{1}{4}\sum_{i}\alpha_i^2\phi(z_i).
$$
Since
$|\psi'(z_i)|^2\hat{g}(\psi(z_i))=e^{\phi(z_i)}\hat{g}(z_i)$ and $\Delta_{\alpha_i}=-\frac{1}{4}\alpha_i\alpha_i+
\frac{Q}{2}\alpha_i$ we conclude
$$
e^{C({\bf \psi(z)})}\prod_i\hat{g}(\psi(z_i))^{\Delta_{\alpha_i}}e^{\frac{s}{4}\sum_i\alpha_i\phi(z_i)}
=e^{C({\bf z})}\prod_i(|\psi'(z_i)|^{-2}\hat{g}(z_i))^{\Delta_{\alpha_i}}.
$$
The proof is completed by the identity
\begin{equation}\label{Dand m}
D_\psi=-\frac{1}{2}m_{\hat{g}}(\phi)
\end{equation}
proven in the appendix. \qed

\subsection{The Liouville measure}\label{sec:LM}

Here, we study the Liouville measure $Z(\cdot )$, the law of which is defined for all Borel sets $A_1, \cdots, A_k\subset \R^2$ by
 \begin{align*}
&  \E^{  \gamma,\mu}_{(z_i,\alpha_i)_i,\hat{g}}  [F(  Z(A_1), \cdots, Z(A_k)      )]  \\
=& (\Pi_{  \gamma,\mu}{(z_i,\alpha_i)_i}(\hat{g},1))^{-1}  \lim_{\epsilon\to 0}\int_\R\E\Big[F( ( e^{\gamma c}\epsilon^{\frac{\gamma^2}{2}}\int_{A_j}e^{\gamma (X_{\hat{g},\epsilon}+ Q/2\ln\hat{g}) } \,d\lambda)_j  )\prod_i \epsilon^{\frac{\alpha_i^2}{2}}e^{\alpha_i (c+X_{\hat{g},\epsilon}+Q/2\ln\hat{g})(z_i)}\nonumber\\
&\exp\Big( - \frac{Q}{4\pi}\int_{\R^2}R_{\hat{g}} (c+X_{\hat{g} } )  \,d\lambda_{\hat{g}} - \mu e^{\gamma c}\epsilon^{\frac{\gamma^2}{2}}\int_{\R^2}e^{\gamma (X_{\hat{g},\epsilon}+ Q/2\ln\hat{g}) }\,d\lambda \Big) \Big]\,dc.\nonumber %
\end{align*}
%
%

%
%
%
%
 
In what follows, we call $Z_0(\cdot)$ the measure   defined under $\P$ by  
$$Z_0(A):=\int_Ae^{\gamma H_{\hat{g}} }\,dM_\gamma $$
so that $Z_0$ in \eqref{Z0def} is $Z_0(\R^2)$. We have:

\begin{proposition}\label{prop:vol}
Under $\P^{  \gamma,\mu}_{(z_i,\alpha_i)_i,\hat{g}}$, the Liouville measure is given for all $A_1, \cdots, A_k$ by    
\begin{equation} \label{defLioumeasure}
\E^{  \gamma,\mu}_{(z_i,\alpha_i)_i,\hat{g}}[F( Z(A_1), \cdots, Z(A_k)    )]=\frac{\int_0^{\infty}  \E\Big[ F\big(  y \frac{Z_0(A_1)}{Z_0(\R^2)}, \cdots,  y \frac{Z_0(A_k)}{Z_0(\R^2)} \big)Z_0(\R^2)^{-\frac{\sum_i \alpha_i -2Q}{\gamma}}  \Big ] e^{- \mu y}  y^{\frac{\sum_i \alpha_i-2Q}{\gamma} -1} dy}{ \mu^{ \frac{2Q-\sum_i \alpha_i}{\gamma} } \Gamma( \frac{\sum_{i} \alpha_i  -2Q}{\gamma}  )\E\Big[  Z_0(\R^2) ^{-\frac{\sum_i\alpha_i-2Q}{\gamma}}\Big]}.
\end{equation}
In particular, \\
1) the volume of the space $Z(\R^2)$ follows the Gamma distribution $\Gamma\Big( \frac{\sum_i \alpha_i-2Q}{\gamma},\mu    \Big) $, meaning 
$$\forall F\in C_b(\R_+),\quad \E^{  \gamma,\mu}_{(z_i,\alpha_i)_i,\hat{g}}[F(Z(\R^2))]=\frac{\mu^{\frac{\sum_i \alpha_i-2Q}{\gamma}}}{\Gamma\Big(\frac{\sum_i \alpha_i -2Q}{\gamma}\Big)}\int_0^{\infty}F(y)y^{\frac{\sum_i \alpha_i -2Q}{\gamma}-1}e^{-\mu y}\,dy.$$
2) the law of the random measure $Z(\cdot)$ conditionally on $Z(\R^2)=A$ is given by
\begin{equation}\label{defUVLM}
\E^{  \gamma,\mu}_{(z_i,\alpha_i)_i,\hat{g}}[F(Z(\cdot))|Z(\R^2)=A)]=\frac{\E\Big[ F\big(A   \frac{Z_0(\cdot)}{Z_0(\R^2)}\big)Z_0(\R^2) ^{-\frac{\sum_i\alpha_i-2Q}{\gamma}}\Big]}{\E\Big[  Z_0(\R^2) ^{-\frac{\sum_i\alpha_i-2Q}{\gamma}}\Big]}
\end{equation}
for any  continuous bounded functional $F$ on the space of finite measures equipped with the topology of weak convergence.\\
3) Under $\P^{  \gamma,\mu}_{(z_i,\alpha_i)_i,\hat{g}}$, the law of the random measure $Z(\cdot)/A$ conditioned on $Z(\R^2)=A$ does  not depend on $A$ and is explicitly given by
$$\E^{  \gamma,\mu}_{(z_i,\alpha_i)_i,\hat{g}}[F(Z(\cdot)/A)|Z(\R^2)=A)]=\frac{\E\Big[ F\big(    \frac{Z_0(\cdot)}{Z_0(\R^2)}\big)Z_0(\R^2) ^{-\frac{\sum_i\alpha_i-2Q}{\gamma}}\Big]}{\E\Big[  Z_0(\R^2) ^{-\frac{\sum_i\alpha_i-2Q}{\gamma}}\Big]}.$$
\end{proposition} 

\noindent {\it Proof.}  Taking the limit $\epsilon\to 0$ in the relation \eqref{eq:insert} gives
\begin{align*}
& \E^{  \gamma,\mu}_{(z_i,\alpha_i)_i,\hat{g}}[F(Z(A_1), \cdots, Z(A_k))]  \\
& = (\Pi^{  \gamma,\mu}_{(z_i,\alpha_i)_i}(\hat{g},1))^{-1}\Big(\prod_i\hat{g}(z_i)^{-\frac{\alpha_i^2}{4}+\frac{Q}{2}\alpha_i}\Big) e^{C(\hat{g})} \\
& \int_\R e^{\big(\sum_i\alpha_i-2Q\big)c}\,\E\Big[F( e^{\gamma c}Z_0(A_1), \cdots , e^{\gamma c}Z_0(A_k)  )  
 \exp\big( - \mu e^{\gamma c}Z_0(\R^2)  \big) \Big]\,dc.
\end{align*}
Finally, let us make the change of variables $e^{\gamma c}Z_0(\R^2)=y $ to complete the proof.\qed

\begin{remark} 
The law of the volume of the sphere $Z(\R^2)$ given above is precisely what one expects to get from scaling limits of planar maps, see Section \ref{sec:maps}.
\end{remark}

%

 \subsection{The unit volume Liouville measure}\label{sec:UVLM}
 
 In the previous section, we introduced the Liouville measure $Z$ and the unit volume Liouville measure; i.e. the law of $Z$ conditionally on $Z(\R^2)=1$ given by \eqref{defUVLM}. These measures exist provided the Seiberg bounds  \eqref{introS} hold. Recall that the first Seiberg bound in \eqref{introS} is a consequence of the integration over the constant mode $c$: this bound entails that the partition function does not diverge when $c \to - \infty$. This condition on the constant mode $c$ disappears when conditioning $Z$ to have unit volume: indeed, this conditioning amounts to fixing the value of $c$. Therefore, it is natural to expect that the unit volume measure can be defined on a larger set of $(z_i,\alpha_i)$; indeed,we have the following:
\begin{lemma}
Suppose that for all $i$, we have $\alpha_i<Q$. The random variable $Z_0(\R^2)$ has a moment of order $\frac{2Q-\sum_i \alpha_i}{\gamma}$ if and only if
\begin{equation}\label{condUVLM}
Q-\frac{\sum_i \alpha_i}{2}< \frac{2}{\gamma}  \wedge \min_i (Q-\alpha_i)
\end{equation}
In particular, under these conditions, the unit volume measure given by \eqref{defUVLM} is well defined.    
\end{lemma}
 
\noindent {\it Proof.} 
 It is standard in the theory of Gaussian multiplicative chaos that $\E[M_\gamma(\R^2)^q] < \infty$ if and only if $q<\frac{4}{\gamma^2}$. Therefore, is is straightforward to see that one must have $Q-\frac{\sum_i \alpha_i}{2}< \frac{2}{\gamma} $. 
 
 Therefore, in the rest of the proof, we will suppose that $\min_i (Q-\alpha_i) < \frac{2}{\gamma} $. To show the if part of the lemma, it suffices to prove that for all $\alpha_j$
 \begin{equation*}
 \E[   (  \int_{B(z_j,1)}  \frac{1}{|x|^{\gamma \alpha_j}} M_\gamma(dx) )^{  \frac{2Q-\sum_i \alpha_i}{\gamma}  }  ]  < \infty.
 \end{equation*}
  By Kahane's convexity inequalities, it is enough to show the above with $M_\gamma$ replaced by a Gaussian multiplicative chaos $\bar{M}_\gamma$ associated to the log-correlated field $X$ with covariance $\E[X(x)X(y)]=\ln_{+} \frac{1}{|x-y|}$. One can construct a cut-off approximation $X_\epsilon$ to $X$ such that for all $\lambda<1$, $(X_{\lambda \epsilon} (\lambda x ))_{|x| \leq 1}= \overset{(Law)}{=} (X_{\epsilon} (x ))_{|x| \leq 1}+\Omega $ where $\Omega$ is an independent centered Gaussian variable with variance $\ln \frac{1}{\lambda}$ (see \cite{review} for instance). Using this cut-off approximation, we have for all $p < \frac{4}{\gamma^2}$
 \begin{align}
 & \E[   (  \int_{B(z_j,1)}  \frac{1}{(|x|+\epsilon)^{\gamma \alpha_j}}   e^{\gamma X_\epsilon(x)- \frac{\gamma^2}{2}  \E[ X_\epsilon(x)^2] })^{ p  }  ]  \nonumber  \\
 &  \leq  \E[   (  \int_{B(z_j,\frac{1}{2})}  \frac{1}{(|x|+\epsilon)^{\gamma \alpha_j}}   e^{\gamma X_\epsilon(x)- \frac{\gamma^2}{2}  \E[ X_\epsilon(x)^2] })^{ p  }  ]   +  \E[   (  \int_{\frac{1}{2} \leq |z_j| \leq 1}  \frac{1}{(|x|+\epsilon)^{\gamma \alpha_j}}   e^{\gamma X_\epsilon(x)- \frac{\gamma^2}{2}  \E[ X_\epsilon(x)^2] })^{ p  }  ]    \nonumber   \\
 & \leq     \E[   (  \int_{B(z_j,\frac{1}{2})}  \frac{1}{(|x|+\epsilon)^{\gamma \alpha_j}}   e^{\gamma X_\epsilon(x)- \frac{\gamma^2}{2}  \E[ X_\epsilon(x)^2] })^{ p  }  ]   + C  \nonumber \\
& \leq  \frac{1}{2^{\xi(p) -\gamma \alpha_j p}} \E[   (  \int_{B(z_j,1)}  \frac{1}{(|x|+2 \epsilon)^{\gamma \alpha_j}}   e^{\gamma X_{2 \epsilon}(x)- \frac{\gamma^2}{2}  \E[ X_ {2 \epsilon}(x)^2] })^{ p  }  ]   + C, \label{ineqrecursiv} 
 \end{align} 
 where $C$ is form line to line a constant independent from everything. Hence, we conclude that 
\begin{equation}\label{ineqproof}
\E \left [   (  \int_{B(z_j,1)}  \frac{1}{(|x|+\epsilon)^{\gamma \alpha_j}}   e^{\gamma X_\epsilon(x)- \frac{\gamma^2}{2}  \E[ X_\epsilon(x)^2] })^{ p  } \right ] dx
\end{equation}  
 is bounded independently of $\epsilon$ if $\xi(p) -\gamma \alpha_j p >0$: this is clear for $\epsilon$ an inverse power of $2$ by applying recursively \eqref{ineqrecursiv}. Otherwise, if $\epsilon$ belongs to a segment $[\frac{1}{2^n},\frac{1}{2^{n-1}} ]$ then from Kahane's convexity inequalities one can bound up to some global multiplicative constant the expectation in \eqref{ineqproof} by the same quantity with $\epsilon$ replaced by $\frac{1}{2^n}$. Now, one can conclude by the fact that $\xi(  \frac{2Q-\sum_i \alpha_i}{\gamma}) -\gamma \alpha_j  \frac{2Q-\sum_i \alpha_i}{\gamma} >0$ is equivalent to $Q-\frac{\sum_i \alpha_i}{2}<  (Q-\alpha_j)$.  This yields one side of the lemma. 
 
 For the only if part of the lemma, along the same lines, one can show that if $ \E[   (  \int_{B(z_j,1)}  \frac{1}{|x|^{\gamma \alpha_i}} \bar{M}_\gamma(dx) )^{ p  }  ]  < \infty$ then we have
 \begin{equation*}
 \E[   (  \int_{B(z_j,1)}  \frac{1}{|x|^{\gamma \alpha_i}} \bar{M}_\gamma(dx) )^{ p  }  ]  \geq  \E[   (  \int_{B(z_j,\frac{1}{2})}  \frac{1}{|x|^{\gamma \alpha_i}} \bar{M}_\gamma(dx) )^{ p  }  ] \geq  \frac{1}{2^{\xi(p) -\gamma \alpha_j p}}   \E[   (  \int_{B(z_j,1)}  \frac{1}{|x|^{\gamma \alpha_i}} \bar{M}_\gamma(dx) )^{ p  }  ] ,
 \end{equation*}    
  hence we get that $\xi(p) -\gamma \alpha_j p >0$.

 \qed

 Therefore, we can define the unit volume Liouviile measure under the condition of the above lemma.

\subsection{Changes of conformal metrics, Weyl anomaly and central charge}\label{sec:Weyl}
In this section, we want to study how  the Liouville partition function \eqref{eq:metricg} depends on the background metric  $g$ conformally equivalent to the spherical metric in the sense of Section \ref{sphermetric}, say $g=e^{\varphi}\hat{g}$. 

\medskip

 By making the change of variables $c \rightarrow c-m_{\hat{g}}(X_{g})$ in \eqref{eq:metricg} and using Proposition \ref{changemetric}, we can and will replace $X_{g}$ by $X_{\hat{g}}$ in the expression \eqref{eq:metricg}.

\medskip
Now we apply the Girsanov transform to the curvature term $e^{- \frac{Q}{4\pi}\int_{\R^2}R_{g} X_{\hat{g} }   \,d\lambda_{g}}$. Since by \eqref{curvature} $R_{g} \lambda_{g}=
(R_{\hat g}-\Delta_{\hat g} \varphi)\lambda_{\hat g}$ this has the effect of shifting the field $X_{\hat{g}}$ by 
$$-\frac{Q}{4\pi}G_{\hat{g}}(R_{\hat g}-\Delta_{\hat g} \varphi)=-\frac{Q}{2}(\varphi-m_{\hat{g}}(\varphi))$$
where we used $G_{\hat{g}}R_{\hat g}=0$ (since $R_{\hat g}$ is constant). 

This Girsanov transform  has also the effect of multiplying the whole partition function by the exponential of 
\begin{align*}
\frac{Q^2}{32\pi^2}&\iint_{\R^2\times\R^2}R_g(z)G_{\hat{g}}(z,z')R_g(z')\,\lambda_g(dz) \lambda_g(dz') \\
=&\frac{Q^2}{16\pi}\int_{\R^2 }R_g(\varphi-m_{\hat{g}}(\varphi))\,d\lambda_g   \\
=&\frac{Q^2}{16\pi}\int_{\R^2 }(R_{\hat{g}}-\triangle_{\hat{g}}\varphi)(\varphi-m_{\hat{g}}(\varphi))\,d\lambda_{\hat{g}}  \quad \text{ (use \eqref{curvature})} \\
=& \frac{Q^2}{16\pi}\int_{\R^2 }| \partial^{\hat{g}} \varphi|^2\,d\lambda_{\hat{g}}    .
\end{align*}
Therefore, by making the change of variables $c \rightarrow c+Q/2 m_{\hat{g}}(\varphi)$ to get rid of the constant $m_{\hat{g}}(\varphi)$ in the expectation, we get
\begin{align}\label{eq:metricg2}
 \Pi_{\gamma,\mu}^{(z_i\alpha_i)_i}  (g,F)=&e^{\frac{1}{96\pi}\int_{\R^2}|\partial^{\hat{g}} \varphi|^2+2R_{\hat{g}}\varphi \,d\lambda_{\hat{g}}+\frac{Q^2}{16\pi}\int_{\R^2 }| \partial^{\hat{g}} \varphi|^2\,d\lambda_{\hat{g}} +Q^2m_{\hat{g}}(\varphi)  } \\
  &\lim_{\epsilon\to 0}\int_\R\E\Big[F( X_{\hat{g} }+c+Q/2\ln \hat{g})\prod_i \epsilon^{\frac{\alpha_i^2}{2}}e^{\alpha_i (c+X_{\hat{g},\epsilon}+Q/2\ln \hat{g})(z_i)}\nonumber\\
&\exp\Big( - \frac{Q}{4\pi}\int_{\R^2}R_{g} c \,d\lambda_{g} - \mu e^{\gamma v}\epsilon^{\frac{\gamma^2}{2}}\int_{\R^2}e^{\gamma X_{\hat{g},\epsilon}+ Q/2\ln \hat{g} }\,d\lambda \Big) \Big]\,dc.\nonumber %
\end{align}
Now we observe that the Gauss-Bonnet theorem  entails
$$\int_{\R^2}R_{g} c \,d\lambda_{g}=\int_{\R^2}R_{\hat{g}} c \,d\lambda_{\hat{g}}$$ because $c$ is a constant. Therefore, using $Q^2m_{\hat{g}}(\varphi)=\frac{6Q^2}{96\pi}\int_{\R^2}2R_{\hat{g}}\varphi\,d\lambda_{\hat{g}} $,
\begin{align}\label{eq:metricg2}
\Pi_{\gamma,\mu}^{(z_i\alpha_i)_i}& (g,F)\\
=&e^{\frac{1+6Q^2}{96\pi}\int_{\R^2}|\partial^{\hat{g}} \varphi|^2+2R_{\hat{g}}\varphi \,d\lambda_{\hat{g}} } \lim_{\epsilon\to 0}\int_\R\E\Big[F( X_{\hat{g} }+c+Q/2\ln \hat{g})\prod_i \epsilon^{\frac{\alpha_i^2}{2}}e^{\alpha_i (c+X_{\hat{g},\epsilon}+Q/2\ln \hat{g})(z_i)}\nonumber\\
&\exp\Big( - \frac{Q}{4\pi}\int_{\R^2}R_{\hat{g}} (c+X_{\hat{g} } ) \,d\lambda_{\hat{g}} - \mu e^{\gamma c}\epsilon^{\frac{\gamma^2}{2}}\int_{\R^2}e^{\gamma X_{\hat{g},\epsilon}+ Q/2\ln \hat{g} }\,d\lambda \Big) \Big]\,dc\nonumber \\
=&e^{\frac{1+6Q^2}{96\pi}\int_{\R^2}|\partial^{\hat{g}} \varphi|^2+2R_{\hat{g}}\varphi \,d\lambda_{\hat{g}} } \Pi_{\gamma,\mu}^{(z_i\alpha_i)_i}  (\hat{g},F).
\end{align}
We can rewrite the above relation in a more classical physics language
\begin{theorem}{\bf (Weyl anomaly and central charge)}
\begin{enumerate}
\item We have the so-called {\bf Weyl anomaly}
$$\Pi_{\gamma,\mu}^{(z_i\alpha_i)_i}  (e^{\varphi}\hat{g},F)=\exp\Big(\frac{c_L}{96\pi}\Big(\int_{\R^2}|\partial \varphi|^2\,d\lambda+ \int_{\R^2}2R_{\hat{g}}\varphi \,d\lambda_{\hat{g}}\Big)\Big) \Pi_{\gamma,\mu}^{(z_i\alpha_i)_i}  (\hat{g},F)$$ where 
$$c_L=1+6Q^2$$ is the {\bf central charge} of the Liouville theory.
\item The law of the Liouville field $\phi$ under $\P^{\gamma,\mu}_{(z_i,\alpha_i)_i,g}$ is independent of the metric $g$ in the conformal equivalence class of $\hat{g}$.
\end{enumerate}
\end{theorem}
Notice that the above theorem can be reformulated as a {\bf Polyakov-Ray-Singer formula } for LQG, see   \cite{R-S} and \cite{Pol,sarnak} for more on this topic.
 
\section{About the $\gamma\geq 2$ branches of Liouville Quantum Gravity} 
Here we discuss various situations that may arise in the study of the case $\gamma\geq 2$.   We want  this discussion to be very concise, so we just give the results as well as references in order to find the tools required  to carry out the computations in full details. Yet, we stress that the computations consist in following verbatim the strategy of this paper.  In what follows, we will only give the partition function in the round metric as the Weyl anomaly then gives straightforwardly the partition function for any metric conformally equivalent to the spherical metric.

\subsection{The case $\gamma=2$ or string theory} 
The case $\gamma=2$ corresponds to $Q=2$ and is very important in string theory, see the excellent review \cite{Kle} as well as the original paper \cite{Pol}. The partition function of LQG is then the limit
\begin{align}\label{eq:defPigcrit}
\Pi_{2,\mu}^{(z_i\alpha_i)_i}& (\hat{g},F)\\
 =& \lim_{\epsilon\to 0}\int_\R\E\Big[F( X_{\hat{g} }+c+\ln\hat{g})\prod_i \epsilon^{\frac{\alpha_i^2}{2}}e^{\alpha_i (c+X_{\hat{g},\epsilon}+ \ln\hat{g})(z_i)}\nonumber\\
&\exp\Big( - \frac{1}{2\pi}\int_{\R^2}R_{\hat{g}} (c+X_{\hat{g},\epsilon})  \,d\lambda_{\hat{g}} - \mu \sqrt{2/\pi} e^{2 c}(-\ln \epsilon)^{1/2}\epsilon^{2}\int_{\R^2}e^{2 X_{\hat{g},\epsilon}+2 \ln\hat{g} }\,d\lambda \Big) \Big]\,dc.\nonumber %
\end{align}
Notice the additional square root $(-\ln \epsilon)^{1/2}$ in order to get a non trivial renormalized interaction term\footnote{The $\sqrt{2/\pi}$ term appears in relation with the results in \cite{Rnew12} to make the $\gamma=2$ case appear as a suitable limit of the $\gamma<2$ case, see Conjecture \ref{conj} below.}. After carrying the same computations than in \eqref{eq:insert} and taking  the limit $\epsilon\to 0$, we get
\begin{align}\label{eq:insertcrit}
\Pi_{2,\mu}^{(z_i\alpha_i)_i} (\hat{g},F)=&\Big(\prod_i\hat{g}(z_i)^{-\frac{\alpha_i^2}{4}+ \alpha_i}\Big) e^{C({\bf z})}\int_\R e^{\big(\sum_i\alpha_i-4\big)c}\,\E\Big[F( c- \theta_{\hat{g}}+X_{\hat{g}} +H_{\hat{g}}+ \ln\hat{g})   \\&\times  
 \exp\Big( - \mu e^{2 c} \int_{\R^2}e^{2 H_{\hat{g}}(x)} \hat{g}(x)M' (dx) \Big) \Big]\,dc,\nonumber %
\end{align}
where the measure $M'(dx)$ is defined by
$$M'(dx)=(2\E[X_{\hat{g}}^2]-X_{\hat{g}})e^{\gamma X_{\hat{g}}-\frac{\gamma^2}{2}\E[X_{\hat{g}}^2]} \,\lambda_{\hat{g}}(dx)$$ and $ C({\bf z})$ defined as in \eqref{def:chatg}.
One can check as in subsection \ref{KPZ} that this partition function is conformally invariant. The convergence of probability of the renormalized measure $(-\ln \epsilon)^{1/2}\epsilon^{2}\int_{\R^2}e^{2 X_{\hat{g},\epsilon}+2 \ln\hat{g} }\,d\lambda$ has been investigated in \cite{Rnew7,Rnew12} when $X_{\hat{g},\epsilon}$ is a white noise decomposition of the field $X_{\hat{g}}$, which can also be taken as a definition of the regularized field. Convergence in law of of the circle average based regularization measure is carried out via the smooth Gaussian approximations introduced in \cite{review}. Establishing the Seiberg bounds needs some extra care and can be handled via the conditioning techniques used in \cite{LBMcrit}.

\subsection{Freezing in LQG} 
For $\gamma>2$ and $Q=2$, one can define
\begin{align}\label{eq:defPigsuper}
\Pi_{\gamma,\mu}^{(z_i\alpha_i)_i}  (\hat{g},F) =&  \lim_{\epsilon\to 0}\int_\R\E\Big[F( X_{\hat{g} }+c+ \ln\hat{g})\prod_i \epsilon^{\frac{\alpha_i^2}{2}}e^{\alpha_i (c+X_{\hat{g},\epsilon}+\ln\hat{g})(z_i)} \\
&\exp\Big( - \frac{1}{2\pi}\int_{\R^2}R_{\hat{g}} (c+X_{\hat{g} } )  \,d\lambda_{\hat{g}} - \mu e^{\gamma c}\epsilon^{2\gamma-2}\int_{\R^2}e^{\gamma X_{\hat{g},\epsilon}+ \gamma \ln\hat{g} }\,d\lambda \Big) \Big]\,dc\nonumber .%
\end{align}
Here we choose to use a white noise regularization of the field $X_{\hat{g}}$ to stick to the framework in \cite{madaule}.
Notice  the unusual power of $\epsilon$ in order to non-trivially renormalize the interaction term, which gets dominated by the near extrema of the 
field $X_{g,\epsilon}$. Under this framework, the convergence in law of the random measures
$$(-\ln\epsilon )^{\frac{3\gamma}{4}}\epsilon^{2\gamma-2}e^{\gamma X_{\hat{g},\epsilon}}dx\to M'_{\frac{2}{\gamma}}(dx)$$
is established in  \cite{madaule}, where $M'_{\frac{2}{\gamma}}(dx)$ is a random measure characterized by
$$\E[e^{M'_{\frac{2}{\gamma}}(f)}] =\E[e^{-c_\gamma\int_{\R^2}f(x)^{\frac{2}{\gamma}}\hat{g}^{-1}(x)M'(dx)}].$$
Hence the convergence in law in the sense of weak convergence of measures
$$(-\ln\epsilon )^{\frac{3\gamma}{4}}\epsilon^{2\gamma-2} e^{\gamma X_{\hat{g},\epsilon}+  \gamma\ln\hat{g} }\,d\lambda\to \hat{g}^\gamma(x)M'_\alpha(dx).$$
We deduce
\begin{align}\label{eq:insert3}
\Pi_{\gamma,\mu}^{(z_i\alpha_i)_i}& (\hat{g},F)\\
 =&\Big(\prod_i\hat{g}(z_i)^{-\frac{\alpha_i^2}{4}+ \alpha_i}\Big) e^{C({\bf z})}\int_\R e^{\big(\sum_i\alpha_i-4\big)c}\,\E\Big[F( c-\frac{\gamma}{2}\theta_{\hat{g}}+X_{\hat{g}} +H_{\hat{g}}+ \ln\hat{g}) \nonumber \\&\times  
 \exp\Big( - \mu e^{\gamma c} \int_{\R^2}e^{\gamma H_{\hat{g}}(x) }\hat{g}(x)\, M'_{\frac{2}{\gamma}}(dx) \Big) \Big]\,dc,\nonumber\\
 =&\Big(\prod_i\hat{g}(z_i)^{-\frac{\alpha_i^2}{4}+ \alpha_i}\Big) e^{C({\bf z})}\int_\R e^{\big(\sum_i\alpha_i-4\big)c}\,\E\Big[F( c-\frac{\gamma}{2}\theta_{\hat{g}}+X_{\hat{g}} +H_{\hat{g}}+ \ln\hat{g}) \nonumber \\&\times  
 \exp\Big( -c_\gamma \mu^{\frac{2}{\gamma}} e^{2 c} \int_{\R^2}e^{2 H_{\hat{g}}(x) }\hat{g}(x)\, M' (dx) \Big) \Big]\,dc,\nonumber %
\end{align}
with $C({\bf z})$ given by  \ref{def:chatg}. Up to the unusual shape of the  cosmological constant, this is exactly the same partition function as in the critical case $\gamma=2$. The difference is here the law of the Liouville measure $M'_{\frac{2}{\gamma}}(dx)$, which can be seen as a $\alpha={\frac{2}{\gamma}}$-stable transform of the derivative martingale $M'$ and is now  purely atomic (see \cite{madaule} for further details).

\subsection{Duality of LQG } 

The basic tools in order to carry out the following computations can be found in \cite{BJRV}.   Define the dual partition function for $\bar{\gamma}>2$ and $Q=\frac{2}{\bar{\gamma}}+\frac{\bar{\gamma}}{2}$ as 
\begin{align}\label{eq:defPigdual}
\bar{\Pi}_{\bar{\gamma},\mu}^{(z_i\alpha_i)_i} (\hat{g},F) 
 =&   \lim_{\epsilon\to 0}\int_\R\E\Big[F( X_{\hat{g} }+c+Q/2\ln\hat{g})\prod_i \epsilon^{\frac{\alpha_i^2}{2}}e^{\alpha_i (c+X_{\hat{g},\epsilon}+Q/2\ln\hat{g})(z_i)}\nonumber\\
&\exp\Big( - \frac{Q}{4\pi}\int_{\R^2}R_{\hat{g}} (c+X_{\hat{g} } )  \,d\lambda_{\hat{g}} - \mu e^{\bar{\gamma} c}\epsilon^{2 }\int_{\R^2}e^{\bar{\gamma} X_{\hat{g},\epsilon}+ \bar{\gamma} Q/2\ln\hat{g} }\,d\lambda_\alpha \Big) \Big]\,dc 
\end{align}
where $\lambda_\alpha$ is a $\alpha$-stable Poisson measure with spatial intensity $\lambda$ and $\alpha=4/\bar{\gamma}^2$. We get \begin{align}\label{eq:insert4}
\bar{\Pi}_{\bar{\gamma},\mu}^{(z_i\alpha_i)_i}& (\hat{g},F)\\
 =&\Big(\prod_i\hat{g}(z_i)^{-\frac{\alpha_i^2}{4}+ \frac{Q}{2}\alpha_i}\Big) e^{C({\bf z})}\int_\R e^{\big(\sum_i\alpha_i-2Q\big)c}\,\E\Big[F( c-\frac{\bar{\gamma}}{2}\theta_{\hat{g}}+X_{\hat{g}} +H_{\hat{g}}+ \frac{Q}{2} \ln\hat{g}) \nonumber \\&\times  
 \exp\Big( - \mu e^{\bar{\gamma} c} \int_{\R^2}e^{\bar{\gamma} H_{\hat{g}}(x) }\hat{g}^{\frac{\bar{\gamma}}{4}}(x)\, S'_\alpha(dx) \Big) \Big]\,dc\nonumber
 \end{align}
with $C({\bf z})$ defined as usual and $S'_\alpha(dx)$ is a stable Poisson random measure with spatial intensity $e^{\gamma X_{g}-\frac{\gamma^2}{2}\E[X_{g}^2]}\,d\lambda$. By computing the expectation we get
\begin{align}\label{eq:insert4}
\bar{\Pi}_{\bar{\gamma},\mu}^{(z_i\alpha_i)_i}  (\hat{g},1) =&\Big(\prod_i\hat{g}(z_i)^{-\frac{\alpha_i^2}{4}+ \frac{Q}{2}\alpha_i}\Big) e^{C({\bf z})}\int_\R e^{\big(\sum_i\alpha_i-2Q\big)c}\,\nonumber \\&\times  \E\Big[ 
 \exp\Big( - \mu^{\frac{\gamma^2}{4}}\frac{4\Gamma(1-\gamma^2/4)}{\gamma^2} e^{ \gamma c} \int_{\R^2}e^{\gamma H_{\hat{g}}(x) }\hat{g}e^{\gamma X_{g}-\frac{\gamma^2}{2}\E[X_{g}^2]}\,d\lambda \Big) \Big]\,dc\nonumber\\
 =&\frac{\mu^{\frac{2Q-\sum_i\alpha_i}{\bar{\gamma}}}}{\mu^{\frac{2Q-\sum_i\alpha_i}{ \gamma}}}\Big(\frac{4\Gamma(1-\gamma^2/4)}{\gamma^2}\Big)^{\frac{2Q-\sum_i\alpha_i}{ \gamma}}\Pi_{\gamma,\mu}^{(z_i\alpha_i)_i}(\hat{g},1).
 \end{align}
Observe that this is an ad-hoc construction of duality (see also \cite{Dup:houches}). The very problem to fully justifies the duality of LQG is to find a proper analytic continuation of the partition of LQG, i.e. the function 
$$\gamma\mapsto  \Pi_{\gamma,\mu}^{(z_i\alpha_i)_i}(\hat{g},1).$$
First observe that this mapping goes to $\infty$ as $\gamma\to 2$ and it is necessary to get rid of the pole at $\gamma=2$. We make the following conjecture
\begin{conjecture}\label{conj}
The function 
$$\gamma\mapsto \Big(\frac{4\Gamma(1-\gamma^2/4)}{\gamma^2}\Big)^{\frac{2Q-\sum_i\alpha_i}{ \gamma}}\Pi_{\gamma,\mu}^{(z_i\alpha_i)_i}(\hat{g},1)$$ is an analytic function of $\gamma\in]0,2[$, which admits an analytic extension for $\gamma\geq 2$ given by $\bar{\Pi}_{\gamma,\mu}^{(z_i\alpha_i)_i}  (\hat{g},1)$. Furthermore, this extension at $\gamma=2$ is the partition function $\Pi_{2,\mu}^{(z_i\alpha_i)_i}(\hat{g},1)$ of the critical case.
\end{conjecture} 
We do not know how to establish analyticity but we stress that the above function is continuous on $]0,+\infty[$.


\section{Perspectives and Conjectures}
In this section, we give a brief overview of perspectives and open problems linked to this work. 

\subsection{The DOZZ formula}\label{DOZZ}

One of the interesting features of LQG is that it is a non minimal CFT but nevertheless physicists have conjectured exact formulas for the three point correlation function of the theory. This correlation function is very important because (in theory) one can compute all correlation functions of LQG from the knowledge of the three point  function.  In LQG, the three point function is quite amazingly supposed to have a completely explicit form, the celebrated DOZZ formula \cite{Do,teschner,ZZ}. 

More precisely, let $z_1, z_2, z_3 \in \R^2$ and $\alpha_1, \alpha_2, \alpha_3$ be three points satisfying the Seiberg bounds \eqref{introS}. Applying the M\"obius transformation rule  \eqref{KPZscaling} for the map
 $\psi $ that  takes $(z_1,z_2,z_3)$ to $(0,1,\infty)$ we get after some
calculation
$$
\Pi_{\gamma,\mu}^{(z_i,\alpha_i)_i} (\hat{g},1)= |z_1-z_2|^{2 \Delta_{12}} |z_2-z_3|^{2 \Delta_{23}} |z_1-z_3|^{2 \Delta_{13}} C_{\gamma}( \alpha_1, \alpha_2, \alpha_3)
$$
where we denoted  $\Delta_{12}= \Delta_{\alpha_3}-\Delta_{\alpha_1}-\Delta_{\alpha_2}  $ and similarly for $\Delta_{13}$ and $\Delta_{23}$. The coefficient is given by (recall  $s=\sum_{i=1}^3 \alpha_i-2Q$)
$$
C_{\gamma}( \alpha_1, \alpha_2, \alpha_3)=e^{\frac{1}{4}(s^2+2Qs)+2\ln 2 \Delta(\alpha_1)} \gamma^{-1}\mu^{-s/\gamma}\Gamma(s/\gamma)
\E\,  Z^{-s/\gamma}
$$
and 
$$
Z=\int |z|^{-\alpha_1\gamma}|z-1|^{-\alpha_2\gamma}\hat g(z)^{-\frac{\gamma}{4}\sum_{i=1}^{3}\alpha_i}M_\gamma(dz)
.$$ 
The DOZZ formula is a conjecture on an exact expression for $C_{\gamma}( \alpha_1, \alpha_2, \alpha_3)$. It is based on the observation that $\E\,  Z^{-s/\gamma}
$ can be computed in closed form if $-s/\gamma=n$, a positive integer. We have
$$
\E\,  Z^n=\int e^{\gamma^2\sum_{i<j}G_{\hat g}(z_i,z_j)} \prod_{i=1}^n|z_i|^{-\alpha_1\gamma}|z_i-1|^{-\alpha_2\gamma}\hat g(z_i)^{-\frac{\gamma}{4}\sum_{j=1}^{3}\alpha_j}\lambda_{\hat g}(dz_i).
$$
Using \eqref{hatGformula} this becomes
$$
\E\,  Z^n=e^{-\gamma^2\frac{n^2-n}{4}}\int \prod_{i<j}|z_i-z_j|^{-\gamma^2} \prod_{i=1}^n|z_i|^{-\alpha_1\gamma}|z_i-1|^{-\alpha_2\gamma}\lambda(dz_i),
$$
an expression that does not depend on the background metric $\hat g$. This Coulomb gas integral can be computed in closed form and  leads to an expression which can be cast in a form where  $n$ enters as a parameter allowing a formal extension of the formula to the negative real axis. This leads to the DOZZ formula for $C_{\gamma}( \alpha_1, \alpha_2, \alpha_3)$. We will not state it here explicitly  as it is quite complicated and involves introducing numerous special functions.

Proving the DOZZ formula seems at this time difficult. Note for instance that for given $\gamma$
only a finite number of positive moments of $Z$ exist so one can  not attempt to solve a moment problem. In the semiclassical $\gamma\to 0$ limit we want to point to an interesting  recent approach to the 
DOZZ formula  by performing deformation of the integration contour in function space \cite{witten}.


\subsection{The semi-classical limit }\label{scl}
The semiclassical limit of LQG is the study of the concentration phenomena of the Liouville field around the extrema of the Liouville action for small $\gamma$, see \cite{nakayama,witten}. After a suitable rescaling of the parameters $\mu$ and $(\alpha_i)_i$, that is 
\begin{equation}\label{asymp}
\mu \gamma^2=\Lambda,\quad \alpha_i=\frac{\chi_i}{\gamma}
\end{equation}
for some fixed constants $\Lambda>0$ and weights $(\chi_i)_i$ satisfying $\chi_i<2$ and $\sum_i\chi_i>4$, the  Liouville field $\gamma \phi$ should converge  in law towards $U+\ln\hat{g}$, where $U$ is the solution of   the classical Liouville equation with sources
\begin{equation}\label{eq:source}
\triangle_{\hat{g}} U -R_{\hat{g}}=2\pi\Lambda e^U- 2\pi \sum_i \chi_i\delta_{z_i},\quad \quad \text{with }\int_{\R^2}e^{U}\,d\lambda_{\hat{g}}= \frac{\sum_i\chi_i-4}{\Lambda},
\end{equation}
hence the name of the theory "Liouville quantum gravity". The reader may consult \cite{lacoin} for some partial results in the "toy model" situation where the zero modes have been turned off.

\subsection{Relation with discretized  $2d$  quantum gravity}\label{sec:maps}

In this Section we will  present some precise conjectures on the connection of our results to the work on discrete models of 2d gravity, randoms surfaces and random maps.

The standard way to discretize 2d quantum gravity coupled to matter fields is to consider a statistical mechanics model (corresponding to a conformal field theory with central charge $c_\textrm{m}$) defined on a  random lattice (or random map), corresponding to the random metric, for instance a random triangulation of the sphere. We formulate below precise mathematical conjectures on the relationship of LQG to that setup.

Let $\mathcal{T}_{N}$ be the set of triangulations of $\S^2$  with $N$ faces  and  $\mathcal{T}_{N,3}$ be the set of triangulations with $N$ faces and $3$ marked faces or points (called roots). 

Next consider a model of statistical physics (matter field) that can be defined  on every   $T\in \mathcal{T}_N$.  The list of such models contains pure gravity (no matter field),  Ising model (a spin $\pm 1$ on each triangle or vertex),  the multicritical discrete spin models (which correspond to the discrete series of the minimal CFT with $1/2\le c_\textrm{m}<1$), the O($N$) dilute and dense loop models with $0\le N<2$, the $q=3$ or $q=4$ Potts models and  discrete models associated to  minimal or rational conformal field theories with central charge $c_\textrm{m}$ such that $-2<c_\textrm{m}\le 1$). We refer to \cite{KosInHandbook} for a review and references.

For $T\in \mathcal{T}_{N}$, define the partition function of the matter field on $T$
$$Z_{\textrm{m}}(T,\beta)=\sum_{C_T}W(C_T,\beta)$$ as a sum of configurations $C_T$ (defined as ensemble of some local or geometric discrete degrees of freedom) over $T$ with positive  local Boltzmann weights $W(C_T,\beta)$. These Boltzmann weights depend on some parameters denoted $\beta$ and these parameters are tuned to their critical point $\beta_c$ such that the statistical model coupled to gravity is critical. At this point, the triangulation $T$ has no marked points. Call $Z_N$ the partition function at criticality for  triangulations of size $N$ 
\begin{equation}
\label{defZN}
Z_N=\sum\limits_{T \in \mathcal{T}_{N,3}} Z_{\textrm{m}}(T,\beta_c),
\end{equation}
where we extend straightforwardly the above definition of $Z_{\textrm{m}}(T,\beta_c)$ to triangulations $T$ with marked points (the marked points play no role in the definition of $Z_{\textrm{m}}(T,\beta_c)$).  
 It is conjectured by physicists that  $Z_N$ diverges as $N$ goes to infinity as (see \cite{Amb})
 \begin{equation}
\label{simZN}
Z_N\sim\ N^{3-(2-\gamma_s)-1}e^{\mu^{\textrm{m}}_c N}(1+o(1))
\end{equation}
with $\mu_c^{\textrm{m}}$ some critical ``cosmological constant'' or ``fugacity'' that depends on the critical model considered, 
 and the string exponent $\gamma_s$ can be explicitly expressed in terms of the central charge $c_\textrm{m}$  of the CFT for the matter field through the relations 
 \begin{equation}
\label{gammaQagain}
2-\gamma_s= \frac{2Q}{\gamma}
\qquad\text{for}
 \quad Q=2/\gamma+\gamma/2=\sqrt{(25-c_\textrm{m})/6}.
\end{equation}

Therefore, for $\bar{\mu}>\mu^{\textrm{m}}_c$, the full partition of the system triangulations+matter field 
 \begin{equation}
\label{Zmub}
Z_{\bar{\mu}}=\sum_Ne^{-\bar{\mu} N}Z_N
\end{equation}
 converges 
and we can sample a random triangulation according to this partition function. We are interested in the regime where the system samples preferably the triangulations with a large number of faces. Notice that for $-2<c_{\textrm{m}} \leq 1$, we have $\sqrt{2}< \gamma \leq 2$ and therefore $-1<\gamma_s\leq 0$. From \eqref{simZN}, we see that  the closer $\bar{\mu}$ is to $\mu^{\textrm{m}}_c$, the larger the typical area of the random triangulation (with 3 marked points) is and for $\mu\sim \mu^{\textrm{m}}_c$, the size of the typical area diverges. Therefore,  we are interested in the limit $\bar{\mu}\to \mu^{\textrm{m}}_c$ in the following regime: we assume that $\bar{\mu}$ depends on a parameter   $a>0$ such that 
\begin{equation}\label{cosmoRPM}
\overline{\mu}=\mu_c^{\textrm{m}}+\mu a^2
\end{equation}
 where $\mu$ is a fixed positive constant.  

Let us now explain how to embed   a triangulation $T\in  \mathcal{T}_{N,3}$ onto the sphere $\S^2$ and define a random measure on $\S^2$ out of it. Following \cite{GilRho} (see also \cite[section 2.2]{Curien}), we can equip such a triangulation with a conformal structure (where each face has the geometry of an equilateral triangle). The uniformization theorem tells us that we can then conformally map the triangulation onto the sphere $\S^2$ and the conformal map is unique if  we pick three distinct points $x_1, x_2 ,x_3$   on the sphere $\S^2$ and demand the map to send the three marked points to $x_1, x_2 ,x_3$.  We denote by $\nu_{T,a}$ the corresponding deterministic measure where each triangle of the sphere is given a volume $a^2$. Concretely, the uniformization provides for each face $t\in T$  a conformal map $\psi_t:  \Delta\to \S^2$ where  $\Delta$ is an equilateral triangle
of volume 1. Then  $\nu_{T,a}(dz)=a^2|(\psi_t^{-1})'(z)|^2dz$ on the image triangle $\psi_t( \Delta)$.  In particular, the volume of the total space $\S^2 $ is $N a^2$. Now, we consider the random measure $\nu_{a,\overline{\mu}}$ defined by 
\begin{equation*}
\E^{a,\overline{\mu}}[  F( \nu_{a,\overline{\mu}} )  ]= \frac{1}{Z_{a}}\sum_N e^{-\overline{\mu} N}\sum_{T \in \mathcal{T}_{N,3}} F(\nu_{T,a} )Z_{\textrm{m}}(T,\beta_c),
\end{equation*} 
 for positive bounded functions $F$ where $Z_a$ is a normalization constant. We denote by $\P^{a,\overline{\mu}}$ the probability law associated to $\E^{a,\overline{\mu}}$.

\medskip
We can now state a precise mathematical conjecture:

\begin{conjecture}\label{conjcartes}
Under $\P^{a,\overline{\mu}}$ and under the relation \eqref{cosmoRPM}, the family of random measures  $(\nu_{a,\overline{\mu}})_{a>0} $ converges in law as $a\to0$ in the space of Radon measures equipped with the topology of weak convergence towards the law of the Liouville measure of  LQG with parameter $\gamma$ given by  \eqref{gammaQagain}, cosmological constant $\mu$ and vertex operators at the points $x_1,x_2,x_3$ with weights $\alpha_i=\gamma$ for all $i$.
\end{conjecture}

Note that  $\nu_{a,\overline{\mu}}(\S^2 ) $ converges in law under $\P^{a,\overline{\mu}}$ as $a\to 0$ towards a $\Gamma(\frac{\sum_i \alpha_i-2Q}{\gamma},\mu)$ distribution with parameter $\gamma$, $\mu$ and $\alpha_i=\gamma$ for all $i$, which corresponds precisely to the law of the volume of the space for  LQG with these parameters (see Subsection \ref{sec:LM}).

\vspace{0.3cm}

\noindent
\textbf{Exemple 1: Pure gravity $c_\textrm{m}=0, \gamma=\sqrt{\frac{8}{3}}$}

Pure gravity corresponds to the case when no matter field is put on the triangulation, in which case $Z_{\textrm{m}}(T,\beta)=\sum_{C_T}W(C_T,\beta)=1$ for all $T$. $Z_N$ thus stands for the cardinal of $\mathcal{T}_{N,3}$ and it is known mathematically since Tutte \cite{Tutte} that 
$$Z_N\sim N^{3-\frac{5}{2}-1}e^{\mu^\textrm{m}_c N}(1+o(1))$$ as $N$ goes to infinity. Notice that $3=\frac{\sum_{i=1}^3 \alpha_i}{\gamma}$ where $\alpha_i=\gamma$ for all $i$   and $\frac{5}{2}= \frac{2Q}{\gamma}$ for $\gamma=\sqrt{\frac{8}{3}}$.

 One can check that $\nu_{a,\overline{\mu}}(\S^2) $ converges in law under $\P^{a,\overline{\mu}}$ as $a\to 0$ towards a $\Gamma(\frac{1}{2},\mu)$ distribution with parameter $\gamma=\sqrt{\frac{8}{3}}$, $\mu$ and $\alpha_i=\gamma$ for all $i$.

\vspace{0.2 cm}

\noindent
\textbf{Exemple 2: Ising model $c_\textrm{m}=\frac{1}{2},\gamma=\sqrt{3}$}

According to the physics literature (see \cite{Amb}),  the partition function of the Ising model on triangulations at criticality   $Z^\textrm{Is}_N$ (corresponding to \eqref{defZN}) should diverge as $N^{3-\frac{7}{3}-1}e^{\mu^{\textrm{Is}}_c N}(1+o(1))$ as $N$ goes to infinity (note that the critical temperature is different on the random lattice models from the regular lattice). Once again,  
notice that $3=\frac{\sum_{i=1}^3 \alpha_i}{\gamma}$ where $\alpha_i=\gamma$ for all $i$  and $\frac{7}{3}= \frac{2Q}{\gamma}$ for $\gamma=\sqrt{3}$.  Again,  $\nu_{a,\overline{\mu}}(\S^2) $ converges in law under $\P^{a,\overline{\mu}}$ as $a\to 0$ towards a $\Gamma(\frac{2}{3},\mu)$ distribution.

\vspace{0.4 cm}

Finally, let us also mention that we could state similar conjectures to conjecture \ref{conjcartes} in the context of fixed volume planar maps. In this context, one samples the map of size $N$ proportionally to the partition function \eqref{defZN} such that it has a fixed volume $A=a^2N$ and then lets $N$ go to infinity (with $a^2= \frac{A}{N}$). The limiting measures will then be (conjecturally) given by the Liouville measure of LQG conditionned to have fixed volume $A$.     

\subsubsection{Conjecture with general vertex operators}
Finally one may ask what is the relation between the general vertex operators $V_\alpha(x)=\exp{(\alpha X(x))}$ (with $\alpha<Q$) that we consider in this paper, the Liouville measure given by \eqref{defLioumeasure} with more than 3 points $x_i$ and some $\alpha_i\neq\gamma$ , and  local observables in discrete 2 dimensional gravity.
Since the 3 original $V_\gamma(x)$ correspond to fixing through conformal invariance the points on $\S^2$, hence to the local density of vertices of the triangulation $T$ through the conformal mapping onto the sphere, it is natural to consider   the local density moment defined as follows.  In addition to the points $x_1,x_2,x_3$ (to which the centers of the marked faces of the triangulation $T$ are sent), we consider additional fixed points $x_i$ with $i>3$ on the sphere, around which  a small disc $\mathcal{D}_{x_i,\epsilon_i}$ centered at $x_i$ with radius $\epsilon_i$ is drawn. Then we   consider the number of vertices $N_{x_i,\epsilon_i}(T)$ of  the triangulation $T$ mapped inside the disk $\mathcal{D}_{x_i,\epsilon_i}$. We consider the random measure defined for all positive bounded functions $F$ as 
\begin{equation*}
\E^{a,\overline{\mu}, (\epsilon_i)_i}[  F( \nu_{a,\overline{\mu},(\epsilon_i)_i} )  ]= \frac{1}{Z_{a,(\epsilon_i)}}\sum_N e^{-\overline{\mu}  N}\sum_{T \in \mathcal{T}_{N,3}} \prod_{i >3}\epsilon_i^{2 \Delta_i} (a^2N_{x_i,\epsilon_i}(T))^{\frac{\alpha_i}{\gamma}} F(\nu_{T,a}) Z_{\textrm{m}}(T,\beta_c),
\end{equation*} 
where $Z_{a,(\epsilon_i)_i}$ is a normalization constant, $\Delta_i=\frac{\alpha_i}{2}(Q-\frac{\alpha_i}{2})$ the conformal weight (see next sections). We denote by $\P^{a,\overline{\mu},(\epsilon_i)_i}$ the probability law associated to $\E^{a,\overline{\mu},(\epsilon_i)_i}$. Notice that we have included the renormalization terms $a^2$ and $\epsilon_i^{2 \Delta_i}$ although  they cancel  with the same
terms in  $Z_{a,(\epsilon_i)_i}$. However, they are needed if one were to consider the limit for the
partition function $Z_{a,(\epsilon_i)_i}$.
We can now state our conjecture:
\medskip

\begin{conjecture}\label{conjcartes2}
Under $\P^{a,\overline{\mu},(\epsilon_i)_i}$ and under the relation \eqref{cosmoRPM}, the family of random measures  $(\nu_{a,\overline{\mu},(\epsilon_i)_i})_{a>0} $ converges in law as $a\to0$ and then as $\epsilon_i \to 0$ in the space of Radon measures equipped with the topology of weak convergence towards the law of the Liouville measure of  LQG with parameter $\gamma$ given by  \eqref{gammaQagain}, cosmological constant $\mu$ and vertex operators at the points $x_1,x_2,x_3$ with weights $\alpha_i=\gamma$ for all $i \leq 3$ and vertex operators at the points $x_i$ with weights $\alpha_i$ for $i>3$.
\end{conjecture}

%
%
%
%
%
\subsubsection{Relation with the Brownian map 
}
It is natural to ask if, in conjecture \ref{conjcartes}, one can reinforce the convergence of measures to a convergence in the space of random metric spaces (equipped with a natural volume form). More precisely, in the case of pure gravity $c_\textrm{m}=0$, 
consider the Riemannian metric defined on each image triangle $\psi_t(\Delta)\subset \S^2$ of the 
uniformization by $a|(\psi_t^{-1})'(z)|^2dz^2$ (hence the lengths of the edges of the image triangles are
$\sqrt{a}$). Let $d_{T,a}$ be the corresponding distance function on $\S_2$ and $d_{a,\bar\mu}$
the random metric on $\S^2$ defined analogously to the random measure $\nu_{a,\bar\mu}$.
Then, it is  widely believed that the metric space (equipped with a volume measure) $(\S^2, d_{a,\bar\mu},\nu_{a,\bar\mu})$ converges in law as $a \to 0$ towards a metric space $(\S^2, d, \nu)$, where $\nu$ is the LQG measure of conjecture \ref{conjcartes}. If this is the case, then the space $(\S^2, d, \nu)$ should be related to the Brownian map equipped with its volume measure (see \cite{LeGall,Mier}): more precisely, for all fixed $A>0$, both metric spaces should be isometric (up to some global constant) once conditioned to have same volume $A$. The isometry should also send the Brownian map volume measure to the measure $\nu$. 



\appendix

\section{M\"obius transform relations}

In this section, we gather a few relations concerning M\"obius transforms and their behavior with respect to Green functions. Recall that the set of automorphisms of the Riemann sphere can be described in terms of the M\"obius transforms
$$\psi(z)=\frac{az+b}{cz+d},\quad a,b,c,d\in \C\text{ and }ad-bc\not=0.$$
Such a function preserves the cross ratios: for all distinct points $z_1,z_2,z_3,z_4\in\C$
\begin{equation}\label{cross}
 \frac{(z_1-z_3)(z_2-z_4)}{(z_2-z_3)(z_1-z_4)}= \frac{(\psi(z_1)-\psi(z_3))(\psi(z_2)-\psi(z_4))}{(\psi(z_2)-\psi(z_3))(\psi(z_1)-\psi(z_4))}.
\end{equation}
Recall that $g_\psi$ stands for the metric $|\psi'|^2   \hat{g}\circ\psi $.
\medskip

\noindent {\it Proof of Proposition \ref{GreenCC}.} We can rewrite the expression \eqref{covXphi} with $g=g_\psi$ in a condensed way
$$G_{g_\psi}(x,y)=\frac{1}{(4\pi)^2}\iint_{\R^2\times\R^2}\ln\frac{|x-z||y-z'|}{|x-y||z-z'|}\,\lambda_{g_\psi}(dz)\lambda_{g_\psi}(dz').$$
By making a change of variables and use \eqref{cross}, we get
\begin{align*}
G_{g_\psi}(x,y)=&\frac{1}{(4\pi)^2}\iint_{\R^2\times\R^2}\ln\frac{|x-\psi^{-1}(z)||y-\psi^{-1}(z')|}{|x-y||\psi^{-1}(z)-\psi^{-1}(z')|}\,\lambda_{\hat{g}}(dz)\lambda_{\hat{g}}(dz').\\
=&\frac{1}{(4\pi)^2}\iint_{\R^2\times\R^2}\ln\frac{|\psi(x)-z||\psi(y)-z'|}{|\psi(x)-\psi(y)||z-z'|}\,\lambda_{\hat{g}}(dz)\lambda_{\hat{g}}(dz').
\end{align*}
This is exactly the expression of $G_{\hat{g}}(\psi(x),\psi(y))$.\qed

\medskip

\begin{corollary}\label{secondterm}
We have the following relations for all M\"obius transforms $\psi$  
\begin{align}
\label{second1}
-2m_{\hat{g}}(\ln\frac{1}{|x-\cdot|})=& -\frac{1}{2}\ln \hat{g}(x)
+\ln 2\\
-2m_{g_\psi}(\ln\frac{1}{|x-\cdot|})+\theta_{g_\psi}
=& -\frac{1}{2}\ln \hat{g}(\psi(x))-\ln|\psi'(x)|+\theta_{\hat{g}}+\ln 2.
\label{second2}
\end{align}
In particular \eqref{hatGformula} holds.
\end{corollary}

\noindent {\it Proof.} We use the following relation
\begin{equation}\label{identite}
\int_{\R^2}   \ln |x-\cdot|   \lambda_{|\psi'|^2   \hat{g}( \psi )} = 2 \pi (\ln ( |ax+b|^2+|cx+d|^2  )-\ln ( |a|^2 +|c|^2)).
\end{equation}
The proof of this identity is based on the fact that both sides have the same Laplacian and the difference of both functions goes to $0$  as $|x|$ goes to infinity.

The first relation is a straightforward consequence of \eqref{identite} with $\psi(z)=z$. One could use  \eqref{identite} as well to prove the second but another way (which we follow below) is to use \eqref{cross}. Write
\begin{align*}
-2m_{g_\psi}(\ln\frac{1}{|x-\cdot|})+\theta_{g_\psi}=&\frac{1}{(4\pi)^2}\int_{\R^2}\int_{\R^2}\ln\frac{|x-z||x-z'|}{|z-z'|}\lambda_{g_\psi}(dz)\lambda_{g_\psi}(dz')\\
=&\frac{1}{(4\pi)^2}\int_{\R^2}\int_{\R^2}\ln\frac{|x-\psi^{-1}(z)||x-\psi^{-1}(z')|}{|\psi^{-1}(z)-\psi^{-1}(z')|}\lambda_{\hat{g}}(dz)\lambda_{\hat{g}}(dz')
\end{align*}
Observe that the mapping $(x,y)\mapsto \frac{1}{(4\pi)^2}\int_{\R^2}\int_{\R^2}\ln\frac{|x-\psi^{-1}(z)||y-\psi^{-1}(z')|}{|\psi^{-1}(z)-\psi^{-1}(z')|}\lambda_{\hat{g}}(dz)\lambda_{\hat{g}}(dz')$ is a continuous function so that we can write
\begin{align*}
-2m_{g_\psi}&(\ln\frac{1}{|x-\cdot|})+\theta_{g_\psi}\\
=&\lim_{y\to x}\frac{1}{(4\pi)^2}\int_{\R^2}\int_{\R^2}\ln\frac{|x-\psi^{-1}(z)||y-\psi^{-1}(z')|}{|\psi^{-1}(z)-\psi^{-1}(z')|}\lambda_{\hat{g}}(dz)\lambda_{\hat{g}}(dz')\\
=&\lim_{y\to x}\big(\frac{1}{(4\pi)^2}\int_{\R^2}\int_{\R^2}\ln\frac{|x-\psi^{-1}(z)||y-\psi^{-1}(z')|}{|x-y||\psi^{-1}(z)-\psi^{-1}(z')|}\lambda_{\hat{g}}(dz)\lambda_{\hat{g}}(dz')+\ln|x-y|\big).
\end{align*}
Now we can use the invariance of cross-products with respect to M\"obius transforms to get
\begin{align*}
-2m_{g_\psi}&(\ln\frac{1}{|x-\cdot|})+\theta_{g_\psi}\\
=&\lim_{y\to x}\big(\frac{1}{(4\pi)^2}\int_{\R^2}\int_{\R^2}\ln\frac{|\psi(x)-z||\psi(y)-z'|}{|\psi(x)-\psi(y)||z-z'|}\lambda_{\hat{g}}(dz)\lambda_{\hat{g}}(dz') +\ln|x-y|\big)\\
=&\lim_{y\to x}\big(\frac{1}{(4\pi)^2}\int_{\R^2}\int_{\R^2}\ln\frac{|\psi(x)-z||\psi(y)-z'|}{|z-z'|}\lambda_{\hat{g}}(dz)\lambda_{\hat{g}}(dz') -\ln\frac{|\psi(x)-\psi(y)|}{|x-y|}\big)\\
&=-2m_{ \hat{g}}(\ln\frac{1}{|\psi(x)-\cdot|})+\theta_{\hat{g}}-\ln|\psi'(x)|.
\end{align*}
We complete the proof thanks to \eqref{second1}.\qed

\begin{lemma}\label{lem:zarb}
The relations \eqref{Grule} and \eqref{Dand m} hold.
\end{lemma}


\proof Using the relation \eqref{identite}, we have
\begin{align*}
&  G_{\hat{g}}(\psi(x),\psi(z))   \\
 =  & \ln \frac{1}{|x-z|}+\frac{1}{2} (\ln ( |ax+b|^2+|cx+d|^2  )-\ln ( |a|^2 +|c|^2))\\
      &  +  \frac{1}{2} (\ln ( |az+b|^2+|cz+d|^2  )-\ln ( |a|^2 +|c|^2))  \\
      & - \frac{1}{4 \pi}  \int_{\R^2}  \frac{1}{2} (\ln ( |au+b|^2+|cu+d|^2  )-\ln ( |a|^2 +|c|^2))    \lambda_{g_\psi} (du)   \\
  =  & \ln \frac{1}{|x-z|}+\frac{1}{2} \ln ( |ax+b|^2+|cx+d|^2  )+  \frac{1}{2} \ln ( |az+b|^2+|cz+d|^2  )-\frac{1}{2}\ln ( |a|^2 +|c|^2)  \\
& - \frac{1}{4 \pi}  \int_{\R^2}  \frac{1}{2} \ln ( |au+b|^2+|cu+d|^2  )     \lambda_{g_\psi} (du) .
\end{align*}
After integrating, we get that
\begin{align*}
   \int_{\R^2}   G_{\hat g}(\psi(x),\psi(z)) \hat{g}(z) dz   =  & -2 \pi \ln (1+ |x|^2)+2 \pi  \ln ( |ax+b|^2+|cx+d|^2  )+\\
 &  \frac{1}{2}  \int_{\R^2}  \ln ( |az+b|^2+|cz+d|^2  )  \lambda_{\hat{g}}(dz)  - 2 \pi \ln ( |a|^2 +|c|^2)  \\
& -   \frac{1}{2} \int_{\R^2}   \ln ( |au+b|^2+|cu+d|^2  )   \lambda_{g_\psi} (du) .
\end{align*}
At this stage, we will suppose that $ad-bc=1$. Hence, we have 
\begin{align*}
     -\frac{1}{2} \int_{\R^2}  & \ln ( |au+b|^2+|cu+d|^2  )   \lambda_{g_\psi} (du)   \\
& =  \frac{1}{2} \int_{\R^2}   \ln (   |\psi'(u)|^2 \hat{g}( \psi(u))  )    \lambda_{g_\psi} (du)    \\
& =  \frac{1}{2} \int_{\R^2}   \ln (   \hat{g}( v)  )    \lambda_{ \hat{g}}( dv)   +  \frac{1}{2} \int_{\R^2}   \ln (   |\psi'(u)|^2  )   \lambda_{g_\psi} (du)    \\
& =  \frac{1}{2} \int_{\R^2}   \ln (   \hat{g}( v)  )      \lambda_{ \hat{g}}( dv)    -   \int_{\R^2}   \ln    |cu+d|    \lambda_{g_\psi} (du) .
\end{align*}
Now, we introduce the function
\begin{equation*}
G(x)= \int_{\R^2}   \ln    |cx+d-cu-d|     \lambda_{g_\psi} (du)  = 4 \pi \ln |c|+\int_{\R^2}   \ln    |x-u|      \lambda_{g_\psi} (du)  .
\end{equation*}
By using equation \eqref{identite}, we get that
\begin{equation*}
G(x)= 4 \pi \ln |c| +  2 \pi (\ln ( |ax+b|^2+|cx+d|^2  )-\ln ( |a|^2 +|c|^2)) 
\end{equation*}
Hence, we get that
\begin{equation*}
\int_{\R^2}   \ln    |cu+d|      \lambda_{g_\psi} (du)  = G(-\frac{d}{c})=4 \pi \ln |c|- 4 \pi \ln |c|-2 \pi  \ln ( |a|^2 +|c|^2)= -2 \pi  \ln ( |a|^2 +|c|^2).
\end{equation*}
At the end, we get
\begin{align*}
\int_{\R^2}  G_{\hat g}(\psi(x),\psi(z)) \lambda_{\hat{g}}(dz)   \\
  = &  -2 \pi \ln (1+ |x|^2)+2 \pi  \ln ( |ax+b|^2+|cx+d|^2  )\\
  &+  \frac{1}{2}  \int_{\R^2}  \ln ( |az+b|^2+|cz+d|^2  )  \lambda_{\hat{g}}(dz) - 2 \pi \ln ( |a|^2 +|c|^2)  \\
& + \frac{1}{2} \int_{\R^2}   \ln (   \hat{g}( v)  )     \lambda_{\hat{g}}(dv) +2 \pi \ln ( |a|^2 +|c|^2)   \\
   = &-  \pi   \ln \frac{g_\psi(x) }{ \hat{g}(x)  } -  \pi m_{\hat g}( \ln \frac{g_\psi(x) }{ \hat{g}(x)  } )  =
   -  \pi\phi(x)-  \pi m_{\hat g}( \phi)
\end{align*}
which implies that
\begin{equation}\label{transs}
\frac{1}{(4\pi)^2} \int_{\R^2}\int_{\R^2}  G_{\hat g}(\psi(x),\psi(z))\lambda_{\hat{g}}(dx)\lambda_{\hat{g}}(dz)=-\frac{1}{2}m_{\hat g}( \phi) . 
\end{equation}
Recall now that $X_{\hat g}\circ\psi$ equals in law $X_{\hat g}-m_{g_\psi}(X_{\hat g})$ so that
$$
 G_{\hat g}(\psi(x),\psi(z))=G_{\hat g}(x,y)-\frac{1}{4\pi}((G_{\hat{g}}e^\phi)(x)+(G_{\hat{g}}e^\phi)(y))+D_\psi
$$
where $
D_\psi=
 \frac{1}{4\pi}m_{\hat g}(e^\phi G_{\hat{g}}e^\phi)$. Using  \eqref{Gonephi} this becomes
\begin{equation}\label{tdpsifinal}
D_\psi=
  \frac{1}{4\pi}( m_{g_\psi}(\phi)-m_{\hat g}(\phi)
)
\end{equation}
and the applying   \eqref{Gonephi} again we get 
$$
 G_{\hat g}(\psi(x),\psi(z))=G_{\hat g}(x,y)-\frac{1}{4}(\phi(x)+\phi(y))+\frac{1}{2}(m_{\hat g}(\phi)+
 m_{g_\psi}(\phi)).
$$
\eqref{transs} implies 
$$
-\frac{1}{2}m_{\hat g}( \phi) =\frac{1}{2}m_{g_\psi}(\phi)
$$
which yields  \eqref{Grule} and combining with \eqref{tdpsifinal} we also get \eqref{Dand m}.
\qed

\end{document}